\documentclass[12pt,a4wide]{amsart}
\usepackage{amsmath,amsthm,amsfonts,amssymb,graphicx}

\usepackage[all]{xy}
\def\demo{{\it Proof. }}
\def\fin{\hfill{$\square$}}
\newtheorem{theorem}{Theorem}[section]
\newtheorem{lemma}[theorem]{Lemma}
\newtheorem{definition}[theorem]{Definition}

\newtheorem{proposition}[theorem]{Proposition}
\newtheorem{corollary}[theorem]{Corollary}
\newtheorem{example}[theorem]{Example}
\newtheorem{algorithm}[theorem]{Algorithm}

\newcommand\Hom{\operatorname{\mathcal{H}om}}
\newcommand\ho{\operatorname{Hom}}

\newcommand\Ker{\operatorname{ker}}

\newcommand\Coh{\operatorname{H}}
\newcommand\coh{\operatorname{h}}
\newcommand\im{\operatorname{Im}}

\newcommand\End{\operatorname{\mathcal{E}nd}}
\newcommand\en{\operatorname{End}}
\newcommand\Spe{\operatorname{Spec}}

\newcommand\ex{\operatorname{Ext}}

\newcommand\Pic{\operatorname{Pic}}

\newcommand\rk{\operatorname{r}}
\newcommand\rg{\operatorname{rk}}
\newcommand\dg{\operatorname{d}}
\newcommand\Jac{\operatorname{\overline{Jac}}}
\newcommand\jac{\operatorname{Jac}}
\newcommand\Ja{\operatorname{\overline{Jac}}}

\newcommand\so{\operatorname{Supp}}
\newcommand\des{\operatorname{\underset{(\leq)}<}}
\newcommand\gen{\operatorname{g}}
\newcommand\Z{\mathbb{Z}}

\title{Simpson Jacobians of reducible curves}
\author{Ana Cristina L\'opez Mart\'{\i}n }
\email{anacris@usal.es}
\address{Departamento de Matem\'aticas, Universidad de Salamanca,
Plaza de la Merced 1-4, 37008 Salamanca, Spain}
\date{\today}
\thanks {e-mail: anacris@usal.es. \\ This research was partly supported by the Spanish-Italian cooperation project
HI00-141, by the research projects BFM2000-1315 and SA009/01 and
by EAGER Contract no. HPRN-CT-2000-00099} 

\begin{document}
\maketitle
\begin{abstract} For any projective curve $X$ let $\overline{M}^d(X)$ be the
Simpson moduli space of pure dimension one rank 1 degree $d$
sheaves that are semistable with respect to a fixed polarization
$H$ on $X$. When $X$ is a reduced curve the connected component of
$\overline{M}^d(X)$ that contains semistable line bundles can be
considered as the compactified Jacobian of $X$. In this paper we
give explicitly the structure of this compactified Simpson
Jacobian for the following projective curves: tree-like curves and
all reduced and reducible curves that can appear as Kodaira
singular fibers of an elliptic fibration, that is, the fibers of
types $III$, $IV$ and $I_N$ with $N\geq 2$
\end{abstract}

\section{Introduction}
The problem of compactifying the (generalized) Jacobian of a
singular curve has been studied since Igusa's work \cite{I} around
1950. He constructed a compactification of the Jacobian of a nodal
and irreducible curve $X$ as the limit of the Jacobians of smooth
curves approaching $X$. Igusa also showed that his
compactification does not depend on the considered family of
smooth curves. An intrinsic characterization of the boundary
points of Igusa's compactification as the torsion free, rank 1
sheaves which are not line bundles is due to Mumford and Mayer.
The complete construction for a family of integral curves over a
noetherian Hensel local ring with residue field separably closed
was carried out by D'Souza \cite{D'So}. One year later, Altman and
Kleiman \cite{AK} gave the construction for a general family of
integral curves.

When the curve $X$ is reducible and nodal, Oda and Seshadri
\cite{OS} produced a family of compactified Jacobians
$\jac_{\phi}$ parameterized by an element $\phi$ of a real vector
space. Seshadri dealt in \cite{Ses1} with the general case of a
reduced curve considering sheaves of higher rank as well.

In 1994, Caporaso showed \cite{C} how to compactify the relative
Jacobian over the moduli of stable curves and described the
boundary points of the compactified Jacobian of a stable curve $X$
as invertible sheaves on certain Deligne-Mumford semistable curves
that have $X$ as a stable model. Recently, Pandharipande \cite{P}
has given  another construction with the boundary points now
representing torsion free, rank 1 sheaves and he showed that
Caporaso's compactification was equivalent to his.

On the other hand, Esteves \cite{E} constructed a compactification
of the relative Jacobian of a family of geometrically reduced and
connected curves and compared it with Seshadri's construction
\cite{Ses1} using theta funtions and Alexeev  \cite{A} gave a
description of the Jacobian of certain singular curves in terms of
the orientations on complete subgraphs of the dual graph of the
curve.

Most of the above papers are devoted to the construction of the
compactified Jacobian of a curve, not to describe it. Moreover
these constructions are only valid for certain projective curves.
However Simpson's work \cite{Si} on the moduli of pure coherent
sheaves on projective schemes allows us to define in a natural way
the Jacobian of any polarized projective curve $X$ as the space
$\jac^d(X)_s$ of equivalence classes of stable invertible sheaves
with degree $d$ with respect to the fixed polarization. This is
precisely the definition we adopt and we also denote by
$\Jac^d(X)$ the space of equivalence classes of semistable pure
dimension one sheaves with rank 1 on every irreducible component
of $X$ and degree $d$.

In some recent papers about the moduli spaces of stable vector
bundles on elliptic fibrations, for instance \cite{HM}, \cite{Br},
\cite{Cal1}, \cite{Cal2}, the Jacobian in the sense of Simpson of
spectral curves appears. Beauville \cite{B} uses it as well in
counting the number of rational curves on K3 surfaces. This
suggests the necessity to determine the structure of these Simpson
Jacobians.

In this paper we give an explicit description of the structure of
these Simpson schemes, $\jac^d(X)_s$ and $\Jac^d(X)$ of any degree
$d$, for $X$ a polarized curve of the following types: tree-like
curves and all reduced and reducible curves that can appear as
singular fibers of an elliptic fibration. The paper is organized
as follows. We will work over an algebraically closed field
$\kappa$ of characteristic zero.

In the second part we define the Simpson Jacobian $\jac^d(X)_s$ of
any projective curve $X$ over $\kappa$, but the definition is also
functorial. By Simpson's work \cite{Si}, the moduli space
$\overline{M}^d(X)$ of equivalence classes of semistable pure
dimension one sheaves on $X$ of (polarized) rank 1 and degree $d$
is a projective scheme that contains $\jac^d(X)_s$ and then it is
a compactification of $\jac^d(X)_s$. The subscheme $\Jac^d(X)$ is
always a connected component of this compactification
$\overline{M}^d(X)$. We show that $\Jac^d(X)$ is only one of the
connected components that can appear in the moduli space
$\overline{M}^d(X)$ (see Proposition \ref{p:otrascomponentes}) and
it coincides with the compactifications constructed by other
authors. The description of the connected components of
$\overline{M}^d(X)$ different to $\Jac^d(X)$ is still an open
problem and, as we show (see Example \ref{e:fraccionario}), it
have to cover the study of moduli spaces of higher rational rank
sheaves on reducible curves. This study is not included in this
paper, however when the curve $X$ is also reduced, as our curves,
the connected component $\Jac^d(X)$ is projective so that it can
be considered as a compactification of the Simpson Jacobian of $X$
as well.

In the third part we collect some general properties about
semistable pure dimension one rank 1 sheaves on reducible curves
and some results relating the Picard groups of two reduced and
projective curves if there is a birational and finite morphism
between them.

In the fourth part we give the description (Theorem
\ref{t:treelike}) of the Simpson Jacobian $\jac^d(X)_s$ and of the
connected component $\Jac^d(X)$ of its compactification when $X$
is a tree-like curve, that is, a projective reduced and connected
curve such that the intersection points of its irreducible
components are disconnecting ordinary double points, but the
singularities lying only at one irreducible component can be
arbitrary singularities. The Picard group of a tree-like curve is
isomorphic to the direct product of the Picard groups of its
irreducible components. Then we use a Teixidor's lemma \cite{T1}
that allows us to order the irreducible components of $X$ and to
find subcurves $X_i$ of $X$ that are also tree-like curves and
intersect their complements in $X$ at just one point. With this
lemma, the determination of the stable rank 1 sheaves on $X$ which
are locally free at all intersection points of the irreducible
components of $X$ is analogous to this given by Teixidor in
\cite{T1} for curves of compact type. When the sheaf is not
locally free at some intersection points, to study its stability
we use a Seshadri's lemma describing the stalk of a pure dimension
one sheaf at the intersection points. Finally for strictly
semistable rank 1 sheaves, we give a recurrent algorithm that
enables us to construct a Jordan-H\"{o}lder filtration and then to
determine their $S$-equivalence classes. The descriptions given
here, and in particular the recurrent algorithm, will be essential
in the next part for the analysis of the Simpson schemes of
Kodaira reduced fibers of an elliptic fibration.

The aim of the fifth part is to describe the structure of the
Simpson Jacobian and of the connected component $\Jac^d(X)$ of its
compactification $\overline{M}^d(X)$ when $X$ is a reducible and
reduced Kodaira fiber of an elliptic fibration, that is, a fiber
of type $III$, $IV$ or $I_N$ with $N\geq 2$. Using the results of
the second section, we first determine the Picard group of these
curves. We find then necessary and sufficient conditions for a
line bundle to be (semi)stable. For the analysis of (semi)stable
pure dimension one rank 1 sheaves which are not line bundles, the
fundamental result is Lemma \ref{l:nodelinea}: if $X$ is a
projective reduced and connected curve, $P$ a singular point at
which $X$ is Gorenstein and $F$ is a pure dimension one rank 1
sheaf on $X$ which is not invertible at $P$, then there is a
projective curve $X'$, a finite birational morphism $\phi\colon
X'\to X$ and a pure dimension one rank 1 sheaf $G$ on $X'$ such
that $\phi_*(G)=F$. Hence we prove that there is an isomorphism
$$\Jac^d(X)-\jac^d(X)\simeq \Jac^{d-1}(X')$$ between the border of the connected component
$\Jac^d(X)$ of $X$ and the whole component $\Jac^{d-1}(X')$ of the
curve $X'$. When $X$ is a fiber of type $III$ or $I_N$, the curve
$X'$ is a tree-like curve and then the structure of this border is
determined by previous results. If $X$ is a fiber of type $IV$,
that is, three rational smooth curves meeting at one point $P$,
the curve $X'$ is also the union of three rational smooth curves
meeting at one point but in such a way that the curve $X'$ cannot
be embedded (even locally) in a smooth surface (see Figure 3). The
study of Simpson schemes for this curve $X'$ is in the subsection
\ref{s:curvaX'}. Lemma \ref{l:localizacion2}, whose proof is a
generalization of the proof given by Seshadri for Lemma
\ref{l:localizacion} in this paper, plays a fundamental role in
the most difficult part of the analysis, the study of (semi)stable
pure dimension one rank 1 sheaves on $X'$ which are not line
bundles. Subsection \ref{s:curvaX'} completes then the description
of the Simpson schemes for the curves of type $IV$.

We finish the paper with a special study for the case of degree
zero. This case is very interesting because the polarization does
not influence  semistability conditions (Corollary
\ref{c:nopolarizacion}). We also prove that if $X$ is a fiber of
type $III$, $IV$ or $I_N$, all strictly semistable pure dimension
one rank 1 sheaves  of degree 0 on $X$ are in the same
$S$-equivalence class (Corollary \ref{c:equivalentes}) and that
the moduli space of semistable pure dimension one sheaves of rank
1 and degree 0 on $X$ is never a fine moduli space (Corollary
\ref{c:nofino}). These results can be found in \cite{Cal1} for a
curve $X$ of type $I_2$ where C\u ald\u araru use them to prove
that the connected component $\Jac^0(X)$ is isomorphic to a
rational curve with one node.

If $p\colon S \to B$ is an elliptic fibration with reduced fibers,
the global structure of the compactified Simpson Jacobian
$\Jac^0(S/B)$ is not totally known (when the fibers of $p$ are
geometrically integral some results can be found in \cite{Br}). We
hope that the descriptions of the compactified Simpson Jacobians
of degree zero given here (Propositions \ref{p:grado0} and
\ref{p:propiasemi}) will be very useful to study the singular
fibers of $\Jac^0(S/B)\to B$, to know the singularities of the
variety $\Jac^0(S/B)$ and to find, using the results of
\cite{Br2}, examples of non isomorphic elliptic fibrations having
isomorphic derived categories.

Two consequences of the given descriptions are the following: the
first is that for reducible curves it is not true in general that
the tensor product of semistable pure dimension one sheaves of
rank 1 is semistable (see Example \ref{e:prodtens}) and the second
is that the pullback of a stable pure dimension one rank 1 sheaf
by a finite morphism of reducible curves is not stable either (see
Example \ref{e:imageninversa}). These two statements constitute
the principal problem to generalize to reducible curves the study
of Abel maps given by Esteves, Gagn\'e and Kleiman in \cite{EGK}
for families of integral curves.

The results presented here are part of my Ph. D. thesis. I am very
grateful to my advisor D. Hern\'andez Ruip\'erez, who introduced
me to the problem of compactifying the generalized Jacobian of a
curve, for his invaluable help and his constant encouragement. I
could not done this work without his support.

\section{The Simpson Jacobian of a projective curve}
Let $X$ be a projective curve over $\kappa$. Let $\mathcal{L}$ be
an ample invertible sheaf on $X$, let $H$ be the associated
polarization and let $h$ be the degree of $H$.

Let $F$ be a coherent sheaf on $X$. We say that $F$ is {\it pure
of dimension one} or {\it torsion free} if for all nonzero
subsheaves $F'\hookrightarrow F$ the dimension of $\so(F')$ is 1.
The (polarized) {\it rank} and {\it degree} with respect to $H$ of
$F$ are the rational numbers $\rk_H(F)$ and $\dg_H(F)$ determined
by the Hilbert polynomial
$$P(F,n,H)=\chi(F\otimes\mathcal{O}_X(nH))=h\rk_H(F)n+\dg_H(F)+\rk_H(F)\chi(\mathcal{O}_X).$$
The {\it slope} of $F$ is defined by
$$\mu_H(F)=\frac{\dg_H(F)}{\rk_H(F)}$$ The sheaf $F$ is {\it
stable} (resp. {\it semistable}) with respect to $H$ if $F$ is
pure of dimension one and for any proper subsheaf
$F'\hookrightarrow F$ one has
$$\mu_H(F')<\mu_H(F) \ (\text{resp.} \leq)$$

In \cite{Si} Simpson defined the multiplicity of $F$ as the
integer number $h\rk_H(F)$ and the slope as the quotient
$$\frac{\dg_H(F)+\rk_H(F)\chi(\mathcal{O}_X)}{h\rk_H(F)}$$
Stability and semistability considered in terms of Simpson's slope
and in terms of $\mu_H$ are equivalent. We adopt these definitions
of rank and degree of $F$ because they coincide with the classical
ones when the curve $X$ is integral. Note however that if $X$ is
not integral, the rank and the degree of a pure dimension one
sheaf are not in general integer numbers.

According to the general theory, for every semistable sheaf $F$
with respect to $H$ there is a {\it Jordan-H\"{o}lder filtration}
$$0=F_0\subset F_1\subset \hdots \subset F_n=F$$ with stable
quotients $F_i/F_{i-1}$ and $\mu_H(F_i/F_{i-1})=\mu_H(F)$ for
$i=1,\hdots,n$. This filtration need not be unique, but {\it the
graded object} $Gr(F)=\textstyle{\bigoplus_{i}} F_i/F_{i-1}$ does
not depend on the choice of the Jordan-H\"{o}lder filtration. Two
semistable sheaves $F$ and $F'$ on $X$ are said to be {\it
$S$-equivalent} if $Gr(F)\simeq Gr(F')$. Observe that two stable
sheaves are $S$-equivalent only if they are isomorphic.

In the relative case, given a scheme $S$ of finite type over
$\kappa$, a projective morphism of schemes $f\colon X\to S$ whose
geometric fibers are curves and a relative polarization $H$, we
define the {\it relative rank} and {\it degree} of a coherent
sheaf $F$ on $X$, flat over $S$, as its rank and degree on fibers,
and we say that $F$ is {\it relatively pure of dimension one}
(resp. {\it stable}, resp. {\it semistable}) if it is flat over
$S$ and if its restriction to every geometric fiber of $f$ is pure
of dimension one (resp. stable, resp. semistable).

Let  ${\bf M}^d(X/S,r)_s$ (resp. $\overline{\bf M}^d(X/S, r)_s$)
be the functor which to any $S$-scheme $T$ associates the set of
equivalence classes of stable locally free (resp. relatively pure
dimension one) sheaves on $X_T=X\underset{S}\times T$ with
relative rank $r$ and degree $d$. Two such sheaves $F$ and $F'$
are said to be equivalent if $F'\simeq F\otimes f_T^*N$, where $N$
is a line bundle on $T$ and $f_T\colon X_T\to T$ is the natural
projection. Similarly, we define the functor ${\bf M}^d(X/S,r)$
(resp. $\overline{\bf M}^d(X/S, r)$) of semistable locally free
(resp. relatively pure dimension one) sheaves of relative rank $r$
and degree $d$.

As a particular case of the Simpson's work \cite{Si}, there exists
a projective scheme $\overline{M}^d(X/S,r)\to S$ which coarsely
represents the functor $\overline{\bf M}^d(X/S, r)$. Rational
points of $\overline{M}^d(X/S,r)$ correspond to $S$-equivalence
classes of semistable torsion free sheaves of rank $r$ and degree
$d$ on $X_s$ ($s\in S$). Moreover, ${\bf M}^d(X/S,r)$ is coarsely
represented by a subscheme $M^d(X/S,r)$ and there are open
subschemes $M^d(X/S,r)_s$ and $\overline{M}^d(X/S,r)_s$ which
represent the other two functors.

\begin{definition} \label{d:jacobiana} {\rm The Simpson Jacobian of degree $d$ of the curve $X$ is
$\jac^d(X)_s=M^d(X/\Spe \kappa,1)_s$. We  denote
$\overline{M}^d(X)=\overline{M}^d(X/\Spe \kappa,1)$. The
projective scheme $\overline{M}^d(X)$ is the compactification of
the Simpson Jacobian of $X$ of degree $d$.}
\end{definition}

When $X$ is an integral curve every torsion free rank 1 sheaf is
stable, and then $\jac^d(X)_s$ is equal to the Picard scheme
$\Pic^d(X)$ and $\overline{M}^d(X)$ coincides with
Altman-Kleiman's compactification \cite{AK}.

In some papers about Jacobians of non irreducible curves, for
instance \cite{Ses1}, \cite{T1}, \cite{T2}, \cite{P}, torsion free
rank 1 sheaves  are considered as those pure dimension one sheaves
having rank 1 on every irreducible component of the curve. If the
notion of rank is given by the Hilbert polynomial, there can be,
depending on the degree of the fixed polarization $H$ on the curve
$X$, pure dimension one sheaves of rank 1 whose restrictions to
some irreducible components of $X$ are not of rank 1.

\begin{example} {\rm Let $X$ be the nodal curve which is a union of two
smooth curves $C_1$ and $C_2$ meeting transversally at one point
$P$. Let $H$ be a polarization on $X$ such that
$\deg(H|_{C_1})=\deg(H|_{C_2})=h$ and let $F$ be a locally free
sheaf on $C_1$ of rank 2 and degree $d$. Let us denote by $i\colon
C_1\hookrightarrow X$ the inclusion map. The sheaf $i_*(F)$ is
pure of dimension one on $X$ and, since $$P(i_*(F), n,
H)=P(F,n,H|_{C_1})=2hn+d+ 2\chi({\mathcal O}_{C_1})\, ,$$ one has
$\rk_H(i_*(F))=1$. However, the restriction of $i_*(F)$ to $C_2$
is a torsion sheaf supported at $P$.}
\end{example}

In order to avoid the confusion, we will say that a sheaf $F$ on
$X$ is of {\it polarized rank 1} if $F$ has rank 1 with the
Hilbert polynomial, that is, $\rk_H(F)=1$, whereas  by {\it rank
1} sheaves we mean those sheaves that have rank 1 on every
irreducible component of $X$.

Semistable pure dimension one rank 1 sheaves of degree $d$ have
polarized rank 1 with respect to any polarization and they are a
connected component of the compactification $\overline{M}^d(X)$ of
the Simpson Jacobian of $X$. Let us denote by $\Jac^d(X)$ this
connected component.

Depending on the degree of the polarization $H$ on $X$,
$\Jac^d(X)$ need not be the unique connected component of the
moduli space $\overline{M}^d(X)$. The following proposition shows
all connected components that can appear in $\overline{M}^d(X)$
when $X$ is a union of two integral curves meeting transversally
only at one point.

If $X$ is a projective and reduced curve with irreducible
components $C_1,\hdots,C_N$ and $P_1,\hdots,P_k$ are  the
intersection points of $C_1,\hdots,C_N$, it is known (see
\cite{Ses1}) that for every pure dimension one sheaf $F$ on $X$
there is an exact sequence $$0\to F \to F_{C_1}\oplus \hdots
\oplus F_{C_N}\to T \to 0$$ where we denote
$F_{C_i}=F|_{C_i}/\text{torsion}$ and $T$ is a torsion sheaf whose
support is contained in the set $\{ P_1,\hdots,P_k\}$.

\begin{proposition} \label{p:otrascomponentes} Let $X=C_1\cup C_2$ be a projective curve with
$C_i$ integral curves for $i=1,2$ and $C_1\cdot C_2=P$. Let $H$ be
a polarization on $X$ of degree $h$ and let $h_{C_i}$ be the
degree of the induced polarization $H_{C_i}$ on $C_i$ for $i=1,2$.
It holds that
\begin{enumerate}\item If $h$ is not a multiple of $h_{C_i}$ for
$i=1,2$, then the only connected component of $\overline{M}^d(X)$
is $\Jac^d(X)$.
\item If $h=rh_{C_i}$ only for one $i=1$ or $2$, then
$$\overline{M}^d(X)=\Jac^d(X)\sqcup \overline{M}^{d_i}(C_i,r)$$
where $d_i=d+\chi(\mathcal{O}_X)-r\chi(\mathcal{O}_{C_i})$.
\item If $h=rh_{C_1}=r'h_{C_2}$, then $$\overline{M}^d(X)=\Jac^d(X)\sqcup\overline{M}^{d_1}(C_1,r)
\sqcup \overline{M}^{d_2}(C_2,r')$$ where
$d_1=d+\chi(\mathcal{O}_X)-r\chi(\mathcal{O}_{C_1})$ and
$d_2=d+\chi(\mathcal{O}_X)-r'\chi(\mathcal{O}_{C_2})$.
\end{enumerate}
\end{proposition}

\begin{proof} Let $F$ be a semistable pure dimension one sheaf on
$X$ of polarized rank 1 and degree $d$ with respect to $H$. From
the above exact sequence, we have
$$h=h_{C_1}\rk_{H_{C_1}}(F_{C_1})+h_{C_2}\rk_{H_{C_2}}(F_{C_2})\,
.$$ Thus, if $\rk_{H_{C_i}}(F_{C_i})>0$ for $i=1,2$, the sheaf $F$
is of rank 1 and its $S$-equivalence class, denoted by $[F]$,
belongs to $\Jac^d(X)$. Since $\rk_{H_{C_i}}(F_{C_i})\in
\mathbb{Z}$, if $\rk_{H_{C_1}}(F_{C_1})$ (resp.
$\rk_{H_{C_2}}(F_{C_2})$) is zero, then $h_{C_2}$ (resp.
$h_{C_1}$) must divide to $h$, namely $h=r'h_{C_2}$ (resp.
$h=rh_{C_1}$). In this case, one has $F\simeq F_{C_2}$ (resp.
$F\simeq F_{C_1}$) so that $F$ is a semistable pure dimension one
sheaf on $C_2$ (resp. $C_1$) of rank $r'$ (resp. $r$) and degree
$d_2$ (resp. $d_1$) with respect to $H_{C_2}$ (resp. $H_{C_1}$).
Therefore, the $S$-equivalence class of $F$ belongs to
$\overline{M}^{d_2}(C_2,r)$ (resp. $\overline{M}^{d_1}(C_1,r)$)
and the result follows.
\end{proof}

The problem of describing the connected components of
$\overline{M}^d(X)$ given by the polarized rank 1 sheaves which
are not of rank 1 have to cover then the analysis of Simpson's
moduli spaces of sheaves of higher rank. Furthermore, when the
number of irreducible components of the curve $X$ is bigger than
two, the compactification $\overline{M}^d(X)$ can even contain
moduli spaces of rational rank sheaves on reducible curves as the
following example shows.

\begin{example}\label{e:fraccionario}{\rm Let $X=C_1\cup C_2\cup C_3$ be a compact type
curve. Let us consider a polarization $H$ on $X$ such that
$h_{C_i}=1$ for $i=1,2,3$. Let $L$ be an invertible sheaf on $C_1$
and let $E$ be a vector bundle of rank 2 on $C_2$. Let us consider
the sheaf $F=i_*(L\oplus E)$ where $i\colon C_1\cup
C_2\hookrightarrow X$ is the inclusion map. The sheaf $F$ is pure
of dimension one and it has polarized rank 1 with respect to $H$.
Moreover, the (semi)stability of $F$ with respect to $H$ is
equivalent to the (semi)stability of $L\oplus E$ with respect to
the polarization  $H_{C_1\cup C_2}$ (see Lemma
\ref{l:estabilidadinducida}). However, $L\oplus E$ is a sheaf of
polarized rank 3/2.}
\end{example}

Nevertheless, when the curve $X$ is projective and reduced, the
(semi)\-stability notion of a sheaf given by the slope $\mu_H$ is
equivalent to the (semi)stability notion given by Seshadri
\cite{Ses1} (it is enough to consider as weights the rational
numbers $a_i=\frac{\deg(H_{C_i})}{\deg(H)}$). Therefore the
connected component $\Jac^d(X)$ of $\overline{M}^d(X)$ coincides
with Seshadri's compactification. Thus, in this case $\Jac^d(X)$
is a projective scheme that contains the Simpson Jacobian of $X$
so that it can be considered as a compactification of
$\jac^d(X)_s$.

On the other hand, if $X$ is a stable curve and
$\overline{P}_{d,X}$ denotes the compactification of the
(generalized) Jacobian constructed by Caporaso in \cite{C}, by
considering as polarization the canonical sheaf of $X$ and by
using essentially the Pandharipande's results of \cite{P}, one
easily proves that there exists a bijective morphism $$\Xi\colon
\overline{P}_{d,X}\to \Jac^d(X)$$ from Caporaso's compactification
to the connected component $\Jac^d(X)$.

Finally, the relation between Oda and Seshadri's compactifications
$\jac_{\phi}(X)$ and the connected component $\Jac^d(X)$, when $X$
is a nodal curve, can be found in \cite{A}.

\section{Torsion free sheaves on reducible curves}

\subsection{General properties of semistable sheaves}

Let $X$ be a projective reduced and connected curve over $\kappa$.
Let $C_1,\hdots,C_N$ denote the irreducible components of $X$ and
$P_1,\hdots,P_k$ the intersection points of $C_1,\hdots,C_N$. Let
$H$ be a polarization on $X$ of degree $h$. Henceforth we shall
use the following notation.

{\bf Notation:} If $F$ is a pure dimension one sheaf on $X$, for
every proper subcurve $D$ of $X$, we will denote by $F_D$ the
restriction of $F$ to $D$ modulo torsion, that is,
$F_D=(F\otimes\mathcal{O}_D)/\text{torsion}$, $\pi_D$ will be the
surjective morphism $F\to F_D$ and $F^D=\Ker\pi_D$. We shall
denote by $h_D$ the degree of the induced polarization $H_D$ on
$D$. If $d=\dg_H(F)$ then we shall write $d_D=\dg_{H_D}(F_D)$. The
complementary subcurve of $D$ in $X$, that is, the closure of
$X-D$, will be denoted by $\overline{D}$. If $g=\gen(X)$ denotes
the arithmetic genus of $X$, that is, the dimension of
$\Coh^1(X,\mathcal{O}_X)$, for any pure dimension one sheaf $F$ on
$X$ of polarized rank 1 and degree $d$ with respect to  $H$, let
$b$, $0\leq b<h$,  be the residue class of $d-g$ modulo $h$ so
that
$$d-g=ht+b.$$ For every proper subcurve $D$ of $X$, we shall write
\begin{equation}\label{notacion}
k_D=\frac{h_D(b+1)}{h}\, .
\end{equation}
If $\beta$ is a real number, we use [$\beta$] to denote the
greatest integer less than or equal to $\beta$.

\bigskip

We collect here some general properties we will use later.

The following lemma, due to Seshadri, describes the stalk of a
pure dimension one sheaf on $X$ at the intersection points $P_i$
that are ordinary double points.

\begin{lemma}\label{l:localizacion}
Let $F$ be a pure dimension one sheaf on $X$. If $P_i$ is an
ordinary double point lying in two irreducible components $C_i^1$
and $C_i^2$, then $$F_{P_i}\simeq \mathcal{O}_{X, P_i}^{a_1}\oplus
\mathcal{O}_{C_i^1, P_i}^{a_2}\oplus \mathcal{O}_{C_i^2,
P_i}^{a_3}$$ where $a_1,\ a_2,\ a_3$ are the integer numbers
determined by:
\begin{align}
&a_1+a_2=\rg(F_{P_i}\underset{\mathcal{O}_{X, P_i}}\otimes
\mathcal{O}_{C^1_i,P_i})\notag\\
&a_1+a_3=\rg(F_{P_i}\underset{\mathcal{O}_{X, P_i}}\otimes
\mathcal{O}_{C^2_i,P_i})\notag\\ &a_1+a_2+a_3=\rg(F_{P_i}\otimes
\kappa)\notag
\end{align}
\end{lemma}

\begin{proof} See \cite{Ses1}, Huiti\^{e}me Partie, Prop. 3.
\end{proof}

\begin{lemma}\label{l:estabilidadinducida} Let $F$ be a pure dimension one sheaf on $X$
supported on a subcurve $D$ of $X$. Then $F$ is stable (resp.
semistable) with respect to $H_D$ if and only if $F$ is stable
(resp. semistable) with respect to $H$.
\end{lemma}

\begin{proof} It follows from the equality $$P(F,n,H)=\chi(i_*F\otimes
\mathcal{O}_X(nH))=\chi(F\otimes \mathcal{O}_D(nH_D))=P(F, n,
H_D)$$ where $i\colon D\hookrightarrow X$ is the inclusion map.
\end{proof}

\begin{lemma} \label{l:buenossubhaces} A torsion free  rank 1 sheaf $F$ on $X$ is stable (resp.
semistable) if and only if $\mu_H(F^D)<\mu_H(F)$ (resp. $\leq$)
for every proper subcurve $D$ of $X$.
\end{lemma}

\begin{proof} Given a subsheaf $G$ of $F$ such that $\so(G)=D\subset X$,
let us consider the complementary subcurve $\overline{D}$ of $D$
in $X$. Since $F_{\overline D}$ is torsion free, we have $G\subset
F^{\overline D}$ with $\rk_H(G)=\rk_H(F^{\overline D})$ so that
$\mu_H(G)\leq\mu_H(F^{\overline D})$ and the result follows.
\end{proof}

\begin{lemma}\label{l:desigualdadesgenerales} Let $L$ be an invertible sheaf on $X$
of degree $d$. Then $L$ is (semi)stable with respect to $H$ if and
only if for every proper connected subcurve $D$ of $X$ the
following inequalities hold:
$$-\chi(\mathcal{O}_D)+h_Dt+k_D\des d_D\des
-\chi(\mathcal{O}_D)+h_Dt+k_D+\alpha_D$$ where $\alpha_D=D\cdot
\overline{D}$ is the intersection multiplicity of $D$ and
$\overline{D}$.
\end{lemma}

\begin{proof} Let us write $d=g+ht+b$. If $L$ is (semi)stable with respect to $H$ and
$D$ is a proper subcurve of $X$, the condition $\mu_H(L^D)\des
\mu_H(L)$ is
$$\frac{hd-hd_D+h_D\chi(\mathcal{O}_X)-\chi(\mathcal{O}_D)}{h-h_D}\des
d$$ which is equivalent to
\begin{equation}\label{des1}
-\chi(\mathcal{O}_D)+h_Dt+k_D\des d_D\, .
\end{equation}
Considering the subsheaf $L^{\overline{D}}$ of $L$, yields
\begin{equation}\label{des2}
-\chi(\mathcal{O}_{\overline{D}})+h_{\overline{D}}t+k_{\overline{D}}\des
d_{\overline{D}}\, .
\end{equation}
Since $X=D\cup \overline{D}$ and $\alpha_D=D\cdot \overline{D}$,
we have $d=d_D+d_{\overline{D}}$, $h=h_D+h_{\overline{D}}$ and
$\chi(\mathcal{O}_X)=\chi(\mathcal{O}_D)+\chi(\mathcal{O}_{\overline{D}})-\alpha_D$.
Then \eqref{des1} and \eqref{des2} give the desired inequalities.

Conversely, if $D$ is a connected subcurve of $X$, by the
left-hand side inequality of the statement, we have
$\mu_H(L^D)\des \mu_H(L)$. Otherwise, this holds for every
connected component of $D$ and it is easy to deduce it for $D$.
The result follows then from the former lemma.
\end{proof}

\subsection{Picard groups and normalizations}
We collect here some results relating the Picard groups of certain
projective reduced curves. They are consequences of the following
proposition due to Grothendieck (\cite{EGA}, Prop. 21.8.5).

\begin{proposition}\label{p:grothendieck} Let $X$ and $X'$ be
projective and reduced curves over $\kappa$. Let $\phi\colon X'\to
X$ be a finite and birational \footnote{By a birational morphism
$X'\to X$ of reducible curves we mean a morphism which is an
isomorphism outside a discrete set of points of $X$} morphism. Let
$U$ be the open subset of $X$ such that $\phi\colon
\phi^{-1}(U)\to U$ is an isomorphism and let $S=X-U$. Let us
denote $ \mathcal{O}'_X=\phi_*( \mathcal{O}_{X'})$. Then, there is
an exact sequence
$$0\to (\prod_{s\in S}\mathcal{O}'^{*}_{X,s}/\mathcal{O}_{X,s}^*)/\im\Coh^0(X',\mathcal{O}^*_{X'})\to
\Pic(X)\xrightarrow{\phi^*}\Pic(X')\to 0\, .$$  If the canonical
morphism $\Coh^0(X,\mathcal{O}_X)\to \Coh^0(X',\mathcal{O}_{X'})$
is bijective, then the kernel of $\phi^*$ is isomorphic to
$\prod_{s\in S}\mathcal{O}'^{*}_{X,s}/\mathcal{O}_{X,s}^*$. \fin
\end{proposition}

\begin{corollary}\label{c:grupoaditivo} Let $X$ and $X'$ be two
projective reduced and connected curves over $\kappa$. Let
$\phi\colon X'\to X$ be a birational morphism which is an
isomorphism outside $P\in X$. If $\phi^{-1}(P)$ is one point $Q\in
X'$ and $\frak{m}_{X',Q}^2\subset \frak{m}_{X,P}$, then the
sequence
$$0\to\mathbb{G}_a\to \Pic(X)\xrightarrow{\phi^*}\Pic(X')\to 0$$
is exact. \fin
\end{corollary}

\begin{corollary}\label{c:odaseshadri}
Let $X=\cup_{i\in I}C_i$ be a projective reduced and connected
curve over $\kappa$. Suppose that the intersection points $\{
P_j\}_{j\in J}$ of its irreducible components are ordinary double
points. Let $X'=\sqcup_{i\in I}C_i$ be the partial normalization
of $X$ at the nodes $\{ P_j\}_{j\in J}$. Then, there is an exact
sequence $$0\to (\kappa^*)^m\to \Pic(X)\to \Pic(X')\to 0$$ where
$m=|J|-|I|+1$. \fin
\end{corollary}

\section{Tree-like curves}
In this section we describe the structure of the Simpson Jacobian
$\jac^d(X)_s$ and of the connected component $\Jac^d(X)$ of its
compactification when $X$ is a tree-like curve. We also give a
recurrent algorithm to determine the $S$-equivalence class of a
semistable sheaf on $X$ which we will use to analyze these Simpson
moduli spaces for reduced fibers of an elliptic fibration.

\begin{definition} {\rm An ordinary double point $P$ of a projective
curve $X$ is a disconnecting point if $X-P$ has two connected
components.}
\end{definition}

\begin{definition}\label{d:tree}{\rm
A tree-like curve is a projective, reduced and connected curve
$X=C_1\cup\hdots\cup C_N$ over $\kappa$ such that the intersection
points, $P_1,\hdots,P_k$, of its irreducible components are
disconnecting ordinary double points.}
\end{definition}

Note that the singularities lying only  at one irreducible
component of a tree-like curve can be arbitrary singularities.

In this section, we will assume that $X$ is a tree-like curve of
arithmetic genus $g$. The number of intersection points of the
irreducible components of $X$ is $k=N-1$ and then, by Corollary
\ref{c:odaseshadri}, one has that
$$\Pic(X)\simeq
\prod_{i=1}^N\Pic(C_i)\, .$$

When the irreducible components of $X$ are smooth, that is, $X$ is
a curve of compact type\footnote{Teixidor's tree-like curves are
better known as compact type curves (see, for instance
\cite{HaM})}, Teixidor (\cite{T1}, Lem.1) proves the following
lemma. Her proof is also valid for any tree-like curve.

\begin{lemma}\label{l:orden}  Let $X=C_1\cup \hdots\cup C_N$ a tree-like curve.
Then, the following statements hold:
\begin{enumerate}
\item It is possible to order the irreducible components of  $X$, so that for every $i\leq N-1$
the subcurve $C_{i+1}\cup \hdots\cup C_N$ is connected.
\item For $i\leq N-1$ there are irreducible components, say $C_{i_1},\hdots,C_{i_k}$, with all
subindices smaller than $i$, such that the subcurve $X_i=C_i\cup
C_{i_1}\cup\hdots \cup C_{i_k}$ is connected and intersects its
complement $\overline{X}_i$ in $X$ in just one point $P_i$.
\end{enumerate} \fin
\end{lemma}

Let $H$ be a polarization on $X$ of degree $h$ and let us suppose
from now on that an ordering of the irreducible components of $X$
as in the lemma has been fixed. Fixing the degree $d$ and using
the notations of \eqref{notacion}, this order allows us to define
inductively integer numbers $d_i^X$ as follows:
\begin{align}
&d_i^X=-\chi(\mathcal{O}_{X_i})+h_{X_i}t+[k_{X_i}]+1-d_{i_1}^X-\hdots-d_{i_k}^X,\
\text{ for }  i=1,\hdots, N-1 \\
&d_N^X=d-d_1^X-\hdots-d_{N-1}^X.\notag
\end{align}

We are now going to modify the above numbers to obtain new numbers
$d_i$ associated with $X$. This is accomplished by a recurrent
algorithm. In order to describe it we start by saying that a
connected subcurve $D=C_{j_1}\cup\hdots\cup C_{j_m}$, $m\geq 1$,
of $X$ ordered according to Lemma \ref{l:orden} is {\it final}
either when the numbers $k_{D_{j_t}}$ are not integers for
$t=1,\hdots, m-1$ or $D$ is irreducible.

If $D$ is a final curve, we define $d_{j_t}$ as follows:

1. if $m>1$, $d_{j_t}=d^D_{j_t}$ for $t=1,\hdots,m$.

2. if $m=1$,
$d_{j_1}=h_{C_{j_1}}t+[k_{C_{j_1}}]-\chi(\mathcal{O}_{C_{j_1}})$.

\begin{algorithm}\label{algoritmo} If the curve $X$ is final, $d_i=d_i^X$
for all $i$. Otherwise, let $i$ be the first index for which
$k_{X_i}\in \Z$. We consider the two connected components, $Y=X_i$
and $\overline{Y}=\overline{X_i}$, of $X-P_i$ and we reorder them
according to Lemma \ref{l:orden}. This induces a new ordering
$P^Y_r$ (resp. $P^{\overline Y}_s$) of the points $P_j$ ($j\ne i$)
in $Y$ (resp. $\overline{Y}$). Then,

a) If $Y$ (resp. $\overline{Y}$) is a final curve, the process
finishes for $Y$ (resp. $\overline{Y}$).

b) If $Y$ is not final, we take the first index $r$ of $Y$ for
which $k_{Y_r}\in \Z$, consider the connected components, $Z$ and
$\overline{Z}$, of $Y-P^Y_r$ and reorder them according to Lemma
\ref{l:orden}. If $Z$ and $\overline{Z}$ are final, the process
finishes for $Y$. Otherwise, we iterate the above argument for
those components that are not final and so on. The process
finishes for $Y$ when all subcurves that we find are final.

c) If $\overline{Y}$ is not final, we take the first index $s$ of
$\overline{Y}$ such that $k_{\overline{Y}_s}\in \Z$, consider the
connected components, $W$ and $\overline{W}$, of
$\overline{Y}-P^{\overline{Y}}_s$ and reorder them according to
Lemma \ref{l:orden}. If $W$ and $\overline{W}$ are final, the
process finishes for $\overline{Y}$. Otherwise, we repeat the
above argument for those components that are not final and so on.
The process finishes for $\overline{Y}$ when all subcurves that we
obtain are final.

The algorithm for $X$ finishes when it finishes for both $Y$ and
$\overline{Y}$.
\end{algorithm}

We can now state the theorem that determines the structure of the
Simpson Jacobian of $X$ and of the schemes $\Jac^d(X)_s$ and
$\Jac^d(X)$.

\begin{theorem} \label{t:treelike} Let $X=C_1\cup\hdots\cup C_N$ , $N\geq 2$, be a tree-like curve.

a) If $k_{X_i}$ is not an integer for every $i\leq N-1$, then
$$\jac^d(X)_s=\prod_{i=1}^{N}\Pic^{d_i^X}(C_i)\ \ \ \text{and}$$
$$\footnote{The inclusion is an equality if all irreducible components of $X$ are smooth}
\jac^d(X)_s\subseteq \Jac^d(X)_s=\Jac^d(X)\simeq \prod_{i=1}^{N}
\Ja^{d_i^X}(C_i)$$ where $d_i^X$ are the above integers.

b) If $k_{X_i}$ is an integer for some $i\leq N-1$, then
$$\jac^d(X)_s=\Jac^d(X)_s=\emptyset\ \ \ \text{and}$$
$$\Jac^d(X)\simeq \prod_{i=1}^{N}\Ja^{d_i}(C_i)$$ where $d_i$ are
the integers constructed with the above algorithm.

\end{theorem}

We postpone the proof of this theorem to subsection
\ref{ss:demostracion} and we now give three examples to illustrate
how Algorithm \ref{algoritmo} works.

\begin{example}\label{e:ejemplo1}{\rm Let $X=C_1\cup\hdots\cup C_N$, $N\geq 2$, be a
tree-like curve with a polarization $H$ whose degree $h$ is a
prime number. Suppose that the irreducible components of $X$ are
ordered according to Lemma \ref{l:orden}. Then, since $h_{X_i}$ is
not divisible by $h$, $k_{X_i}=\frac{h_{X_i}(b+1)}{h}$ is an
integer if and only if $b=h-1$. Therefore, by Theorem
\ref{t:treelike}, we have that if $b<h-1$,
$\jac^d(X)_s=\prod_{i=1}^{N}\Pic^{d_i^X}(C_i)$ and
$$\jac^d(X)_s\subsetneqq \Jac^d(X)_s=\Jac^d(X)\simeq
\prod_{i=1}^{N} \Ja^{d_i^X}(C_i),$$ whereas for $b=h-1$,
$\jac^d(X)_s$ and $\Jac^d(X)_s$ are empty and $$\Jac^d(X)\simeq
\prod_{i=1}^{N}\Ja^{d_i}(C_i).$$ Moreover, in this case the number
$d_i$ is given by $$d_i=h_{C_i}(t+1)-\chi(\mathcal{O}_{C_i}) \ \ \
\text{ for } i=1,\hdots,N.$$ Actually, if $L$ is a strictly
semistable line bundle of degree $d$ on $X$, since $k_{D}$ is
integer for all $D\subset X$, the first index $i$ such that
$k_{X_i}$ is integer is $i=1$ and the connected components of
$X-P_1$ are $Y=C_1$ and $\overline{Y}=C_2\cup \hdots \cup C_N$.
Since $C_1$ is irreducible, it is a final subcurve so that
$d_1=h_{C_1}(t+1)-\chi(\mathcal{O}_{C_1})$.  Since
$k_{\overline{Y}_s}$ is integer for $s=2,\hdots,N-1$,
$\overline{Y}$ is not a final curve so we have to apply the
procedure to the curve $\overline{Y}$. Here, the first index $s$
such that $k_{\overline{Y}_s}$ is integer is $s=2$ and the
connected components of $\overline{Y}-P_2$ are $W=C_2$ and
$\overline{W}=C_3\cup\hdots\cup C_N$. Then $d_2=
h_{C_2}(t+1)-\chi(\mathcal{O}_{C_2})$. Since $k_{\overline{W}_s}$
is integer for $s=3,\hdots,N-1$, we apply the procedure to the
curve $\overline{W}$, and so on. The iteration of this procedure
will only finish when we obtain a sheaf supported on $C_N$ that
belongs to $\Pic^{h_{C_N}(t+1)-\chi(\mathcal{O}_{C_N})}(C_N)$.
Thus, $d_i=h_{C_i}(t+1)-\chi(\mathcal{O}_{C_i})$ for
$i=1,\hdots,N$.}
\end{example}

\begin{example}{\rm We are now going to recover examples 1 and 2 of \cite{A}.
There, the irreducible components $C_i$ are taken to be smooth and
$d=g-1$. Then, the residue class of $d-g$ modulo $h$ is $b=h-1$
and $t=-1$. It follows that $k_{X_i}$ is integer for
$i=1,\hdots,N$ and, arguing as in the former example, we have that
$\jac^{g-1}(X)_s$ and $\Jac^{g-1}(X)_s$ are empty and
$$\Jac^{g-1}(X)\simeq
\prod_{i=1}^N\Jac^{h_{C_i}(t+1)-\chi(\mathcal{O}_{C_i})}(C_i)\simeq
\prod_{i=1}^N\Pic^{g_i-1}(C_i)$$ as asserted in \cite{A}.}
\end{example}

\begin{example} {\rm Let $X$ be the following tree-like curve

\vspace{1truecm}
\hspace{4truecm}\includegraphics[scale=0.5]{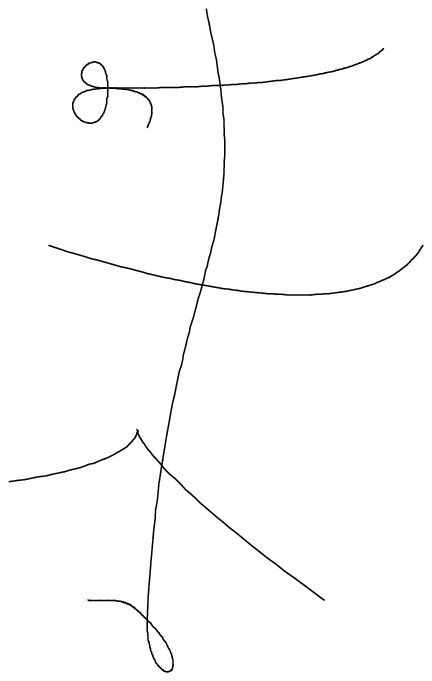}
\vspace{1truecm}

Fixing an ordering of the irreducible components of $X$ as in
Lemma \ref{l:orden}, we obtain

\vspace{1truecm}
\hspace{4truecm}\includegraphics[scale=0.5]{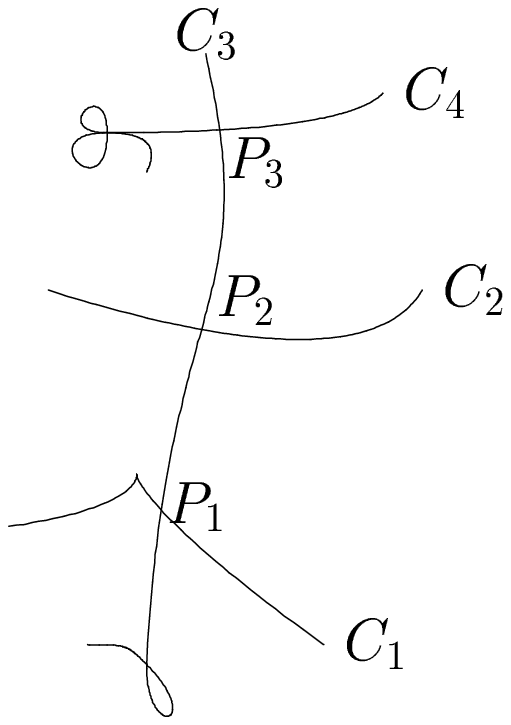}
\vspace{1truecm}

\noindent so that $X_1=C_1$, $X_2=C_2$, $X_3=C_1\cup C_2\cup C_3$.
Assume that the first index $i$ such that $k_{X_i}$ is integer is
$i=3$ and let us compute the numbers $d_i$, $i=1,\hdots,4$ in this
case.

The connected components of $X-P_3$ are $Y=C_1\cup C_2\cup C_3$
and $\overline{Y}=C_4$. Since $C_4$ is a final subcurve, we have
that $d_4=h_{C_4}t+k_{C_4}-\chi(\mathcal{O}_{C_4})$.  We now have
to fix a new ordering for $Y$ according to Lemma \ref{l:orden}. We
can take, for instance, $Y=C_{\sigma(1)}\cup C_{\sigma(2)}\cup
C_{\sigma(3)}$ with $\sigma(1)=1$, $\sigma(2)=3$ and
$\sigma(3)=2$. Then, $Y_{\sigma(1)}=C_1$ and
$Y_{\sigma(2)}=C_1\cup C_3$. Therefore, since
$k_{Y_{\sigma(1)}}=k_{X_1}$ and $k_{Y_{\sigma(2)}}=k_Y-k_{X_2}$
are not integers, $Y$ is final and we conclude that
\begin{align}
d_1&=d^Y_{\sigma(1)}=-\chi(\mathcal{O}_{Y_{\sigma(1)}})+h_{Y_{\sigma(1)}}t+
[k_{Y_{\sigma(1)}}]+1=d^X_1,\notag\\
d_3&=d_{\sigma(2)}^Y=-\chi(\mathcal{O}_{Y_{\sigma(2)}})+h_{Y_{\sigma(2)}}t+
[k_{Y_{\sigma(2)}}]+1-d^Y_{\sigma(1)}=\notag\\
&=-\chi(\mathcal{O}_Y)+h_Yt+k_Y+1-d^X_1-d^X_2-1=d^X_3-1,\notag\\
d_2&=d_{\sigma(3)}^Y=d_Y-d_{\sigma(1)}^Y-d_{\sigma(2)}^Y=d^X_2.\notag
\end{align}}
\end{example}

\subsection{Proof of Theorem
\ref{t:treelike}}\label{ss:demostracion}

In order to prove that the Simpson Jacobian of $X$ is not empty
only when $k_{X_i}\notin \Z$ for all $i\leq N-1$, we need the next
two lemmas that characterize stable invertible sheaves on $X$.
These two lemmas are essentially  steps 1 and 2 in \cite{T1}.

\begin{lemma}\label{l:lineatree} Let $L$ be a line bundle on $X$ of degree
$d$. If $L$ is stable, then $k_{X_i}$ is not an integer for every
$i\leq N-1$ and $L$ is obtained by gluing invertible sheaves $L_i$
on $C_i$ of degrees $d_i^X$, $i=1,\hdots,N$.
\end{lemma}

\begin{proof} Considering the subcurves $X_i$ of $X$,
$i=1,\hdots,N-1$, given by Lemma \ref{l:orden}, from Lemma
\ref{l:desigualdadesgenerales}, we get
\begin{equation}\label{desigualdad}
-\chi({\mathcal O}_{X_i})+h_{X_i}t+k_{X_i}\des d_{X_i}\des
-\chi({\mathcal O}_{X_i})+h_{X_i}t+k_{X_i}+1
\end{equation}
because $\alpha_{X_i}=1$. We have $\rk_{H_{X_i}}(L_{X_i})=1$ so
that $d_{X_i}=\dg_{H_{X_i}}(L_{X_i})$ is an integer. Then, if
$k_{X_i}\in \Z$ for some $i\leq N-1$, \eqref{desigualdad} becomes
a contradiction. Thus $k_{X_i}\notin \Z$ for all $i\leq N-1$ and
there is only one possibility for $d_{X_i}$, namely
$$d_{X_i}=-\chi(\mathcal{O}_{X_i})+h_{X_i}t+[k_{X_i}]+1,\ \ \
\text{ for } i=1,\hdots, N-1. $$ From
$d_{X_i}=d_{C_i}+d_{C_{i_1}}+\hdots+d_{C_{i_k}}$ and the exact
sequence $$0\to L\to L_{C_1}\oplus\hdots\oplus L_{C_N}\to
\oplus_{i=1}^{N-1}\kappa(P_i)\to 0$$ we deduce that
$d_{C_i}=d_i^X$ for all $i$ and the proof is complete.
\end{proof}

The following lemma  proves the converse of Lemma
\ref{l:lineatree}.

\begin{lemma}\label{l:reciproco} Let $L$ be the invertible sheaf on $X$ of degree
$d=g+ht+b$ obtained by gluing line bundles $L_i$ on $C_i$ of
degrees $d_i^X$, $i=1,\hdots,N$. Suppose that $k_{X_i}$ is not
integer for every $i\leq N-1$. Then $L$ is stable.
\end{lemma}

\begin{proof} Taking as weights $\frac{\deg(H_{C_i})}{\deg(H)}$
and proceeding as in Lemma 2 of \cite{T1}, we obtain that for
every proper connected subcurve $D$ of $X$ the following
inequalities are satisfied:
$$-\chi(\mathcal{O}_D)+h_Dt+k_D< d_D<
-\chi(\mathcal{O}_D)+h_Dt+k_D+\alpha_D$$ where $\alpha_D=D\cdot
\overline{D}$. The stability of $L$ follows then from Lemma
\ref{l:desigualdadesgenerales}.

\end{proof}

By Lemmas \ref{l:lineatree} and \ref{l:reciproco}, $\jac^d(X)_s$
is not empty only if $k_{X_i}$ is not integer for every $i\leq
N-1$, and in this case it is equal to
$\prod_{i=1}^N\Pic^{d_i^X}(C_i)$.

We prove now the remaining statements of the theorem. If $L$ is a
strictly semistable line bundle on $X$ of degree $d$ then,
$d_{X_i}$ is equal to one of the two extremal values of the
inequality \eqref{desigualdad}. In particular, $k_{X_i}$ is an
integer for some $i\leq N-1$.

Let $i$ be the first index such that $k_{X_i}$ is integer. Then,
there are two possibilities for $d_{X_i}$:

a) $d_{X_i}=-\chi(\mathcal{O}_{X_i})+h_{X_i}t+k_{X_i}$

b) $d_{X_i}=-\chi(\mathcal{O}_{X_i})+h_{X_i}t+k_{X_i}+1$

Let us construct a Jordan-H\"{o}lder filtration for $L$ in both cases.
Since case a) and case b) are the same but with the roles of $X_i$
and $\overline{X_i}$ intertwined, we give the construction in the
case a).

We have that $\mu_H(L_{X_i})=\mu_H(L^{X_i})=\mu_H(L)$. Then,
$L_{X_i}$ and $L^{X_i}\simeq L_{\overline{X_i}}(-P_i)$ are
semistable with respect to $H$ and, by Lemma
\ref{l:estabilidadinducida}, they are semistable with respect to
$H_{X_i}$ and $H_{\overline{X_i}}$ respectively.

For simplicity, we shall write  $Y=X_i=C_{i_0}\cup
C_{i_1}\cup\hdots \cup C_{i_k}$ with $i_1,\hdots, i_k<i_0=i$ and
$Z=\overline{X_i}$, which are again  tree-like curves.

Let us see when the sheaves $L_Y$ and $L_Z(-P_i)$ are stable. We
can fix an ordering for $Y$ as in Lemma \ref{l:orden}, so that
$Y=C_{\sigma(i_0)}\cup\hdots\cup C_{\sigma(i_k)}$ and we obtain
subcurves  $Y_r$ of $Y$ for
$r=\sigma(i_0),\hdots,\sigma(i_{k-1})$.

{\it Claim 1.} The sheaf $L_Y$ is stable if and only if $k_{Y_r}$
is not an integer for $r=\sigma(i_0),\hdots,\sigma(i_{k-1})$.

\demo Since the residue class of $d_Y-\gen(Y)$ modulo $h_Y$ is
$b_Y=k_Y-1$, the numbers $\frac{h_{Y_r}(b_Y+1)}{h_Y}=k_{Y_r}$ are
not integers for $r=\sigma(i_0),\hdots,\sigma(i_{k-1})$. Then,
from Lemma \ref{l:reciproco}, we have only to prove that $L_Y$ is
in $\prod_r \Pic^{d_r^Y}(C_r)$, where $d_r^Y$ are the integer
numbers defined as $d_i^X$ but with the new ordering of $Y$ and
$r$ runs through the irreducible components of $Y$. This is
equivalent to proving that
$$d_{Y_r}=-\chi(\mathcal{O}_{Y_r})+h_{Y_r}t+[k_{Y_r}]+1\ \ \
\text{ for } r=\sigma(i_0),\hdots,\sigma(i_{k-1}).$$ Actually,
since $L$ is semistable and $Y_r$ is a proper subcurve of $X$, by
Lemma \ref{l:desigualdadesgenerales}, we obtain
\begin{equation}\label{des7}
-\chi(\mathcal{O}_{Y_r})+h_{Y_r}t+k_{Y_r}\leq d_{Y_r}\leq
-\chi(\mathcal{O}_{Y_r})+h_{Y_r}t+k_{Y_r}+\alpha\tag{7}
\end{equation}
where $\alpha$ is the number of intersection points of $Y_r$ and
its complement in $X$. We have that $\alpha\leq 2$ and $d_{Y_r}$
is not equal to the extremal values of \eqref{des7} because
$k_{Y_r}\notin \Z$. Moreover, if it were
$$d_{Y_r}=-\chi(\mathcal{O}_{Y_r})+h_{Y_r}t+[k_{Y_r}]+2,$$ since
$d_Y=d_{Y_r}+d_{\overline{Y_r}^Y}$,
$h_Y=h_{Y_r}+h_{\overline{Y_r}^Y}$ and
$\chi(\mathcal{O}_Y)=\chi(\mathcal{O}_{Y_r})+\chi(\mathcal{O}_{\overline{Y_r}^Y})-1$,
$\overline{Y_r}^Y$ being the complement of $Y_r$ in $Y$, then
$$d_{\overline{Y_r}^Y}=-\chi(\mathcal{O}_{\overline{Y_r}^Y})+h_{\overline{Y_r}^Y}t+
[k_{\overline{Y_r}^Y}]$$ which contradicts  the semistability of
$L$. Thus,
$$d_{Y_r}=-\chi(\mathcal{O}_{Y_r})+h_{Y_r}t+[k_{Y_r}]+1$$ and the
proof of the claim 1 is complete.

On the other hand, the irreducible components of $Z$ are ordered
according the instructions in Lemma \ref{l:orden} and the
subcurves $Z_s$, where $s$ runs through the irreducible components
of $Z$ and $s\leq N-1$, are equal to either $X_s$ or $X_s-Y$.

{\it Claim 2.} The sheaf $L_Z(-P_i)$ is stable if and only if
$k_{X_s}$ is not an integer for every $s>i$.

\demo Since the residue class of $d_Z-1-\gen(Z)$ modulo $h_Z$ is
$b_Z=k_Z-1$ and $k_Y\in \Z$ , by the hypothesis, the numbers
$\frac{h_{Z_s}(b_Z+1)}{h_Z}$ are not integers for $s>i $ and, by
the choice of $i$, they aren't for $s<i$ either. Then, by Lemma
\ref{l:reciproco}, it is enough to prove that
$$d_{H_{Z_s}}(L_Z(-P_i)|_{Z_s})=-\chi(\mathcal{O}_{Z_s})+h_{Z_s}t+[k_{Z_s}]+1\
\ \ \text{ for } s\leq N-1.$$ Since $L$ is semistable, we have
that $$d_{X_s}=-\chi(\mathcal{O}_{X_s})+h_{X_s}t+[k_{X_s}]+1.$$
Moreover, if  $Z_s=X_s$ then,
$d_{H_{Z_s}}(L_Z(-P_i)|_{Z_s})=d_{X_s}$ and if $Z_s=X_s-Y$ then,
$d_{H_{Z_s}}(L_Z(-P_i)|{Z_s})=d_{X_s}-d_Y-1$. We obtain the
desired result in both cases and the proof of the claim 2 is
complete.

We return now to the proof of the theorem. If $k_{Y_r}$ and
$k_{X_s}$ are not integers for
$r=\sigma(i_0),\hdots,\sigma(i_{k-1})$ and $s>i$ (i.e. $Y$ and $Z$
are final curves), then $0\subset L_Z(-P_i)\subset L$  is a
Jordan-H\"{o}lder filtration for $L$ and the $S$-equivalence class of
$L$ belongs to $\prod_r\Pic^{d_r^Y}(C_r)\times
\prod_s\Pic^{d_s^X}(C_s)$.

On the other hand, if $k_{Y_r}$ is integer for some
$r=\sigma(i_0),\hdots,\sigma(i_{k-1})$, the sheaf $L_Y$ is
strictly semistable and we have to repeat the above procedure with
$L_Y$ in the place of $L$ and the curve $Y$ in the place of $X$.
Similarly, if $k_{X_s}$ is integer for some $s>i$, the sheaf
$L_Z(-P_i)$ is strictly semistable. Then, we have to repeat the
above procedure for $L_Z(-P_i)$. By iterating this procedure, we
get a Jordan-H\"{o}lder filtration for $L_Y$: $$0=F_0\subset
F_1\subset \hdots \subset F_m=L_Y$$ and another for $L_Z(-P_i)$:
$$0=G_0\subset G_1\subset \hdots\subset G_n=L_Z(-P_i).$$
Therefore, a filtration for $L$ is given by $$0=G_0\subset
G_1\subset \hdots\subset L_Z(-P_i)\subset \pi_Y^{-1}(F_1)\subset
\hdots\subset \pi_Y^{-1}(L_Y)=L.$$ Thus, the $S$-equivalence class
of $L$ belongs to $\prod_{i=1}^{N}\Pic^{d_i}(C_i)$, where $d_i$
are the integer numbers constructed with the algorithm.

Finally, let us consider a pure dimension one sheaf $F$ on $X$ of
rank 1 and degree $d$ which is not a line bundle. When $F$ is
locally free at the intersection points $P_i$ for all
$i=1,\hdots,N-1$, calculations and results are analogous to the
former ones. If $F$ is not locally free at $P_i$ for some
$i=1,\hdots,N-1$, then there is a natural morphism $$F\to
F_Y\oplus F_Z$$ where $Y$, $Z$ are the connected components of
$X-{P_i}$, that is clearly an isomorphism outside $P_i$. But this
is an isomorphism at $P_i$ as well because by Lemma
\ref{l:localizacion}, $F_{P_i}\simeq \mathcal{O}_{C_i^1,
P_i}\oplus \mathcal{O}_{C_i^2, P_i}$ and this is precisely the
stalk of $F_Y\oplus F_Z$ at $P_i$. We conclude that $F$ is not
stable and if it is semistable, then $k_Y$ and $k_Z$ are integers,
$d_Y$ and $d_Z$ are given by
$$d_Y=-\chi(\mathcal{O}_Y)+h_Yt+k_Y,\quad\quad
d_Z=-\chi(\mathcal{O}_Z)+h_Zt+k_Z,$$ and $F_Y$ and $F_Z$  are
semistable as well. Then, the construction of  a Jordan-H\"{o}lder
filtration for $F$ can be done as above and thus the
$S$-equivalence class of $F$ belongs to
$\prod_{i=1}^{N}\Jac^{d_i}(C_i)$. \fin

\

We give here an example that shows that for reducible curves the
tensor product of two semistable torsion free sheaves of rank 1 is
not in general a semistable sheaf, even when the considered
sheaves are line bundles.

\begin{example}\label{e:prodtens} {\rm Let $X=C_1\cup C_2$ be a tree-like curve with
$C_1\cdot C_2=P$. Suppose that $C_i$ are rational curves and fix
on $X$ a polarization $H$ such that $P\notin \so(H)$ and
$h_{C_i}=1$ for $i=1,2$. Let $Q_1$ and $Q_2$ be two smooth points
of $C_1$. From Lemma \ref{l:desigualdadesgenerales}, it is easy to
prove that the sheaves $\mathcal{O}_X(Q_1)$ and
$\mathcal{O}_X(Q_2)$ are semistable with respect to $H$. The same
lemma proves that $\mathcal{O}_X(Q_1)\otimes \mathcal{O}_X(Q_2)$
is not a semistable sheaf.}
\end{example}

\section{Kodaira reduced fibers} Let $B$ be a projective smooth curve over $\kappa$
and let $p\colon S\to B$ be an elliptic fibration. By this we mean
a proper flat morphism of schemes whose fibers are Gorenstein
curves of arithmetic genus 1. By a Kodaira's result (Thm. 6.2
\cite{Ko}) the singular reduced fibers of $p$ can be classified as
follows:

$I_1: X=C_1$ a rational curve with one node.

$I_2: X=C_1\cup C_2$, where $C_1$ and $C_2$ are rational smooth
curves with $C_1\cdot C_2=P+Q$.

$I_N: X=C_1\cup C_2\cup \hdots \cup C_N$,  $N=3,4,\hdots$, where
$C_i$, $i=1,\hdots, N$, are rational smooth curves and $C_1\cdot
C_2=C_2\cdot C_3=\hdots=C_{N-1}\cdot C_N=C_N\cdot C_1=1$.

$II: X=C_1$ a rational curve with one cusp.

$III: X=C_1\cup C_2$ where $C_1$ and $C_2$ are rational smooth
curves with $C_1\cdot C_2=2P$.

$IV: X=C_1\cup C_2\cup C_3$, where $C_1, C_2, C_3$ are rational
smooth curves and $C_1\cdot C_2=C_2\cdot C_3=C_3\cdot C_1=P$.

In this section we give the description of the Simpson Jacobian
and of the connected component $\Jac^d(X)$ of its compactification
$\overline{M}^d(X)$ for all reduced singular fibers of an elliptic
fibration. Since when the fiber $X$ is an irreducible curve (i.e.
a rational curve with one node or one cusp) every pure dimension
one rank 1 sheaf on $X$ is stable, we mean the fibers of types
$III$, $IV$ and $I_N$ for $N\geq 2$.

\subsection{Preliminary results}
Lemmas \ref{l:morfismo} and \ref{l:gorenstein} can be found  in
\cite{EGK} for integral curves and the proofs given are valid for
finite morphisms of reducible curves.

\begin{lemma}\label{l:morfismo}
Let $\phi\colon X'\to X$ be a finite birational morphism of
projective reduced and connected curves over $\kappa$. Let $G_1$
and $G_2$ be two  pure dimension one rank 1 sheaves on $X'$ and
let $u\colon \phi_*(G_1)\to \phi_*(G_2)$ be a morphism. Then there
is a unique morphism  $v\colon G_1\to G_2$ such that
$\phi_*(v)=u$. \fin
\end{lemma}

\begin{lemma}\label{l:gorenstein} Let $X$ be a projective reduced and connected curve over
$\kappa$ and $P$ a singular point of $X$. Let $\frak{m}$ denote
the maximal ideal of $P$. Set
$\frak{B}:=\End_{\mathcal{O}_X}(\frak{m})$ and
$X^*:=\Spe(\frak{B})$. Let $\psi\colon X^*\to X$ denote the
natural map. Then the following assertions hold:
\begin{enumerate}
\item $\frak{B}=\Hom_{\mathcal{O}_X}(\frak{m},\mathcal{O}_X)$.
\item The curve $X$ is Gorenstein at $P$ if and only if $\psi\colon X^*\to X$
is a finite birational morphism and $\gen(X^*)=\gen(X)-1$.
\item If  $X$ is Gorenstein at $P$ and $\phi\colon X'\to X$
is a birational morphism nontrivial at $P$, then $\phi$ factors
trough $\psi$.
\end{enumerate}\fin
\end{lemma}

Moreover, we can adapt the proof of Lemma 3.8 in \cite{EGK} to
show the following:

\begin{lemma}\label{l:nodelinea}  Let $X$ be a projective reduced and connected curve over $\kappa$
and $P$ a singular point at which $X$ is Gorenstein. Let
$\frak{m}$ denote the maximal ideal of $P$. Set
$\frak{B}:=\End_{\mathcal{O}_X}(\frak{m})$ and
$X^*:=\Spe(\frak{B})$. Let $\psi\colon X^*\to X$ denote the
natural map. Let $F$ be a pure dimension one rank 1 sheaf on $X$.
Then, $F$ is not locally free at $P$ if and only if there is a
pure dimension one rank 1 sheaf $G$ on $X^*$ such that
$\psi_*(G)=F$. If $G$ exists, then it is unique.
\end{lemma}

\begin{proof} If $G$ exits, it is unique by Lemma
\ref{l:morfismo}. The sheaf $G$ exits if and only if $F$ is a
$\frak{B}$-module, so if and only if $\End_{\mathcal{O}_X}(F)$
contains $\frak{B}$.

To reduce the notation, we set $\mathcal{O}:=\mathcal{O}_{X,P}$,
$M:=\frak{m}_P$, $B:=\frak{B}_P$ and $F:=F_P$.

Hence $G$ does not exist if $F$ is invertible at $P$ because then
$\en_\mathcal{O} (F)=\mathcal{O}$, whereas the cokernel of
$\mathcal{O} \hookrightarrow B$ has, by (2) of Lemma
\ref{l:gorenstein}, length 1.

Suppose now that $F$ is not invertible. We have to prove that
$BF\subset F$. Set $F^*:=\ho_\mathcal{O}(F,\mathcal{O})$. Since
the curve $X$ is Gorenstein at $P$, then $F^{**}=F$ (see, for
example, \cite{Ei}). Let $\tilde{X}$ be the total normalization of
$X$, that is, if $X=\cup_{i}C_i$, then
$\tilde{X}=\sqcup_{i}\tilde{C_i}$ where $\tilde{C_i}$ is the
normalization of the integral curve $C_i$. Let us denote
$\overline{\mathcal{O}}:=\mathcal{O}_{\tilde{X},P}$ and
$K:=\overline{\mathcal{O}}^*$.

Since the $\overline{\mathcal{O}}$-module
$$F^*\overline{\mathcal{O}}:=(F^*\otimes_{\mathcal{O}}\overline{\mathcal{
O}})/\text{torsion}$$ is free of rank 1, there is an element $g\in
F^*$ such that
$F^*\overline{\mathcal{O}}=g\overline{\mathcal{O}}$. Then we have
$g\mathcal{O}\subset F^*\subset g\overline{\mathcal{O}}$. By
applying $\ho_{\mathcal{O}}(-,g\mathcal{O})$, we get
\begin{equation}\label{e1}
K\subset gF\subset \mathcal{O}\, .
\end{equation}
Since $F$ is not invertible,  $gF\neq \mathcal{O}$. Hence
$gF\subset M$.

Let $\overline{M}$ be the Jacobson radical of
$\overline{\mathcal{O}}$ which is  invertible as
$\overline{\mathcal{O}}$-module. Set
$G:=(\overline{M}^{-1}\otimes_{\overline{\mathcal{O}}}K)/\text{torsion}$.
The natural morphism
$$\overline{M}^{-1}\otimes_{\overline{\mathcal{O}}}K\otimes_{\overline{\mathcal{O}}}gF\to K$$
induces an inclusion
$$gGF=(G\otimes_{\overline{\mathcal{O}}}gF)/\text{torsion}\hookrightarrow K\, .$$
So \eqref{e1} yields $gGF\subset gF$. Therefore,
\begin{equation}\label{e2}
GF\subset F\, .
\end{equation}
Note that $G\nsubseteq\mathcal{O}$. Indeed, otherwise $G\subset K$
because $K$ is the largest $\overline{\mathcal{O}}$-submodule of
$\mathcal{O}$. However, if $G\hookrightarrow K$, there is a
morphism
$$\overline{M}^{-1}\otimes_{\overline{\mathcal{O}}}K\to K$$ and
then
$$K\hookrightarrow
\overline{M}K:=(\overline{M}\otimes_{\overline{\mathcal{O}}}K)\,
.$$ Hence, by Nakayama's Lemma, $K=0$ which is not true as Lemma
\ref{l:Knonulo}, that we will see later, proves.

The natural morphism
$$\overline{M}^{-1}\otimes_{\overline{\mathcal{O}}}K\to
\ho_{\overline{\mathcal{O}}}(\overline{M},\ho_{\mathcal{O}}(\overline{\mathcal{O}},{\mathcal{O}}))\simeq\ho_{\mathcal{O}}
(\overline{M},\mathcal{O})\, ,$$ gives an inclusion
$$G\hookrightarrow \ho_{\mathcal{O}}(\overline{M},\mathcal{O})/\text{torsion}\, .$$
Since $\mathcal{O}\hookrightarrow
\ho_{\mathcal{O}}(\overline{M},\mathcal{O})/\text{torsion}$, it is
possible to consider $D:=G+\mathcal{O}$. Then $\mathcal{O}\subset
D$, but $\mathcal{O}\neq D$. Set
$D^*:=\ho_{\mathcal{O}}(D,\mathcal{O})$. Since $X$ is Gorenstein,
$D^{**}=D$ and we have $D^*\subset \mathcal{O}$ and $D^*\neq
\mathcal{O}$. Hence $D^*\subset M$. Therefore
$\ho_{\mathcal{O}}(M,\mathcal{O})\subset D$. Since, by (1) of
Lemma \ref{l:gorenstein}, it is
$\ho_{\mathcal{O}}(M,\mathcal{O})=B$, we have
$$BF\subset DF=F+GF\, .$$ Hence $\eqref{e2}$ implies
$BF\subset F$ and the proof is complete.
\end{proof}

\begin{lemma}\label{l:Knonulo}
Using the above notations, one has $K\neq 0$.
\end{lemma}

\begin{proof} Since $\tilde{X}=\sqcup_{i}\tilde{C_i}$, we have $\overline{\mathcal{O}}=\oplus_i{\mathcal{O}}_i$
with ${\mathcal{O}}_i=\mathcal{O}_{{\tilde{C_i}},P}$. Then
$K=\overline{\mathcal{O}}^*=\oplus_i{\mathcal{O}}_i^*$ and it is
enough to prove that ${\mathcal{O}}_i^*\neq 0$.

Using the local duality, one has
$$\ex^1_{\mathcal{O}}({\mathcal{O}}_i, \mathcal{O})^*\simeq
\ho(\coh^0_P({\mathcal{O}}_i),\coh^1_P(\mathcal{O}))\, .$$ Since
${\mathcal{O}}_i$ has depth one, $\coh^0_P(\mathcal{O}_i)=0$ and
then $\ex^1_{\mathcal{O}}(\mathcal{O}_i, \mathcal{O})=0$.

The exact sequence
$$0\to \mathcal{O}\to \oplus_i\mathcal{O}_i\to T\to 0$$ yields
$$0\to \ho_{\mathcal{O}}(\mathcal{O}_i,\mathcal{O})\to \mathcal{O}_i\to
\ho_{\mathcal{O}}(\mathcal{O}_i,T)\to 0$$ because
$\ho_{\mathcal{O}}(\mathcal{O}_i,\mathcal{O}_j)=0$ for $i\neq j$.

Since $T$ is a torsion sheaf, we conclude that
$\ho_{\mathcal{O}}(\mathcal{O}_i,\mathcal{O})\neq 0$ and the
result follows.
\end{proof}

\begin{lemma}\label{l:estporimdi} Let $\phi\colon X'\to X$
be a finite birational and surjective morphism of projective
reduced and connected curves with $\gen(X')=\gen(X)-s$. Let $H$ be
a polarization on $X$ of degree $h$ such that $H':=\phi^*(H)$ has
degree $h$. Let $G$ be a pure dimension one sheaf on $X'$ of rank
1 and degree $d$ with respect to $H'$. Then the following
statements hold:
\begin{enumerate}
\item The sheaf $\phi_*(G)$ is pure of dimension one of rank 1 and degree $d+s$ with respect to $H$.
\item $G$ is (semi)stable with respect to $H'$ if and only if $\phi_*(G)$ es (semi)stable with respect to $H$.
\end{enumerate}
\end{lemma}

\begin{proof} The sheaf $\phi_*(G)$ is pure of dimension one
because it is the direct image of a torsion free sheaf by a finite
morphism. Moreover, for any sheaf $G'$ on $X'$, we have
$P(\phi_*(G'),n,H)=P(G',n,H')$, so that
\begin{align}
\rk_{H}(\phi_*(G'))&=\rk_{H'}(G')\, ,\notag\\
\dg_{H}(\phi_*(G'))&=\dg_{H'}(G')+\rk_{H'}(G')(\chi(\mathcal{O}_{X'})-\chi(\mathcal{O}_X))\,
.\notag
\end{align}
Then
$$\mu_{H}(\phi_*(G'))=\mu_{H'}(G')+s$$ and the result is now straightforward.
\end{proof}

\subsection{The description for the fibers of type $III$}
Here $X$ will denote a fiber of type $III$ of an elliptic
fibration, that is, $X=C_1\cup C_2$ with $C_1$ and $C_2$ two
rational smooth curves and $C_1\cdot C_2=2P$ (figure 1).

\vspace{0.5truecm}
\hspace{4.5truecm}\includegraphics[scale=0.3]{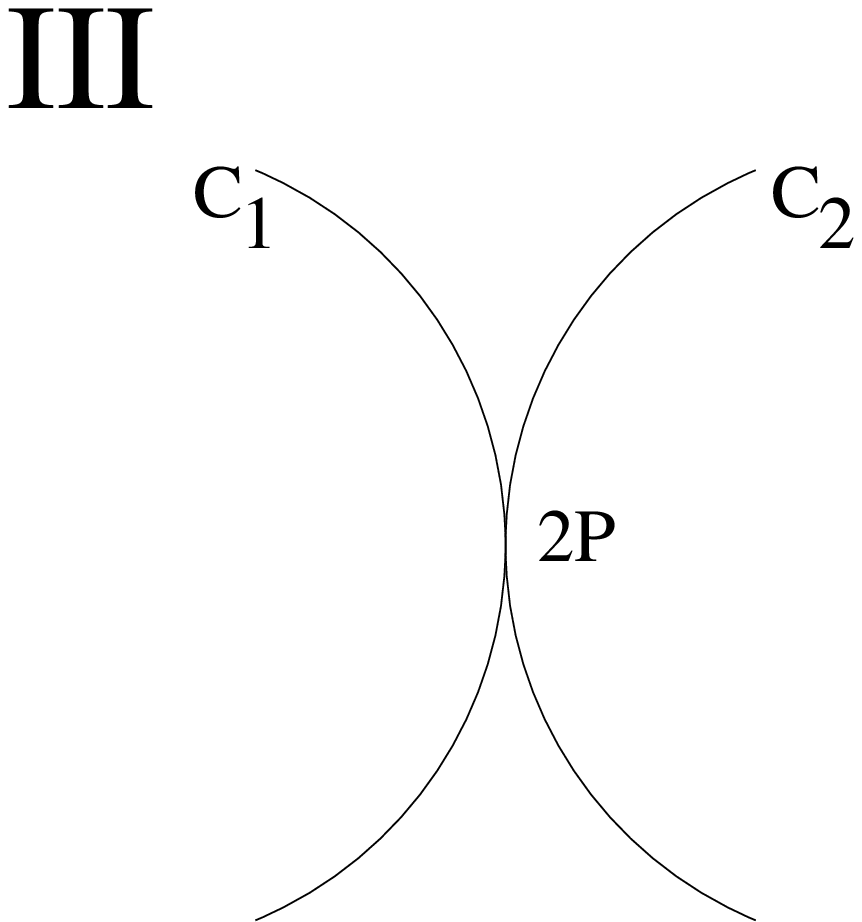}
\vspace{0.5truecm}

\centerline{Figure 1.} \vspace{0.5truecm}

If $X'=C_1\cup C_2$ is the blow-up of $X$ at the point $2P$, then
$X'$ is a tree-like curve with $C_1\cdot C_2=P$ and there is a
finite birational morphism $\phi\colon X'\to X$ which is an
isomorphism outside $2P$ such that $\phi^{-1}(2P)=P$. Since it is
possible to write the completions of the local rings as
$$\widehat{{\mathcal{O}}}_{X,2P}=\kappa[[x,y]]/x(x-y^2)\quad
\text{ and }\quad
\widehat{{\mathcal{O}}}_{X',P}=\kappa[[y,\lambda]]/\lambda(\lambda-y)$$
with $x=\lambda y$, we have that $\frak{m}_{X',P}^2\subset
\frak{m}_{X,2P}$. Hence, by Corollary \ref{c:grupoaditivo} and
taking into account that $X'$ is a tree-like curve, there is an
exact sequence
\begin{equation} \label{e:picIII} 0\to \mathbb{G}_a\to \Pic(X)\to \prod
_{i=1}^2\Pic(C_i)\to 0
\end{equation}
where the last morphism is given by $L\to (L_{C_1}, L_{C_2})$.

Let $H$ be a polarization on $X$ of degree $h$. Fixing the degree
$d$ and using the notations of \eqref{notacion} (in this case
$g=1$), the proposition describing the structure of the Simpson
Jacobian of $X$ and of the scheme $\jac^d(X)$ of $S$-equivalence
classes of semistable line bundles is the following:

\begin{proposition}\label{l:lineaIII}
Let $X=C_1\cup C_2$ be a curve of type $III$.

a) If $k_{C_1}\in \Z$, then there is an exact sequence
$$0\to \mathbb{G}_a\to \jac^d(X)_s\to\prod_{i=1}^2\Pic^{h_{C_i}t+k_{C_i}}(C_i)\to 0\,
, \text{ and }$$
$$\jac^d(X)-\jac^d(X)_s=\prod_{i=1}^2\Pic^{h_{C_i}t+k_{C_i}-1}(C_i)\, .$$

 b) If $k_{C_1}\notin \Z$, then there is an exact sequence
$$0\to \mathbb{G}_a\to \jac^d(X)_s\to\bigsqcup_{i,j=1}^2\Pic^{h_{C_i}t+[k_{C_i}]}(C_i)\times
\Pic^{h_{C_j}t+[k_{C_j}]+1}(C_j)\to 0\, .$$ In this case,
$$\jac^d(X)-\jac^d(X)_s=\emptyset\, ,$$ that is, there are not strictly semistable line bundles.
\end{proposition}

\begin{proof} Since
$\chi(\mathcal{O}_X)=\chi(\mathcal{O}_{C_1})+\chi(\mathcal{O}_{C_2})-2$,
if $L$ is a line bundle on $X$ of degree $d$, one has the exact
sequence $$0\to L\to L_{C_1}\oplus L_{C_2}\to T\to 0$$ where $T$
is a torsion sheaf with support at $2P$ and $\chi(T)=2$.

By Lemma \ref{l:desigualdadesgenerales}, the sheaf $L$ is
(semi)stable if and only if
\begin{equation}\label{d:d1}
h_{C_1}t+k_{C_1}-1\des d_{C_1} \des h_{C_1}t+k_{C_1}+1\, .
\end{equation}
Since $d_{C_1}\in \Z$ and $d=d_{C_1}+d_{C_2}$, we get the
following:

a) if $k_{C_1}\in \Z$, then $L$ is stable if and only if
$d_{C_i}=h_{C_i}t+k_{C_i}$ for $i=1,2$.

b) if $k_{C_1}\notin \Z$, then $L$ is stable if and only if
$d_{C_1}=h_{C_1}t+[k_{C_1}]$ and $d_{C_2}=h_{C_2}t+[k_{C_2}]+1$ or
the same but with the roles of $C_1$ and $C_2$ intertwined.

From \eqref{e:picIII}, we obtain then the exact sequences in a)
and b).

Suppose now that the sheaf $L$ is strictly semistable. In this
case, $d_{C_1}$ is equal to one of the two extremal values of
\eqref{d:d1}. Assume that $d_{C_1}=h_{C_1}t+k_{C_1}-1$ (the other
case is the same but with the roles of $C_1$ and $C_2$
intertwined). Then $k_{C_1}\in \Z$,
$\mu_{H}(L_{C_1})=\mu_{H}(L^{C_1})=d$ and $L_{C_1}$ and $L^{C_1}$
are stable sheaves with respect to $H_{C_1}$ and $H_{C_2}$
respectively. Thus, $L^{C_1}\subset L$ is a Jordan-H\"{o}lder
filtration for $L$ and then its $S$-equivalence class belongs to
$\prod_{i=1}^2\Pic^{h_{C_i}t+k_{C_i}-1}(C_i)$. Hence, there is
only one $S$-equivalence class of strictly semistable line bundles
so that
$\jac^d(X)-\jac^d(X)_s=\prod_{i=1}^2\Pic^{h_{C_i}t+k_{C_i}-1}(C_i)$.
\end{proof}

As before, let $X'=C_1\cup C_2$ be the tree-like curve obtained by
gluing transversally at $P$ the irreducible components of the
curve $X$ and $\phi\colon X'\to X$ the natural morphism. Let $H$
be a polarization on $X$ of degree $h$ such that $H':=\phi^*(H)$
is of degree $h$ as well.

The following proposition proves that there is an isomorphism
between the set of boundary points of the connected component
$\Jac^d(X)$ and the component $\Jac^{d-1}(X')$ when we consider
the polarization $H$ on $X$ and the polarization $H'$ on $X'$.
Since $X'$ is a tree-like curve, the structure of the border
$\Jac^d(X)-\jac^d(X)$ is determined by Theorem \ref{t:treelike}.

\begin{proposition}\label{p:bordeIII} Let $X=C_1\cup C_2$ be a curve of type $III$ with $C_1\cdot C_2=2P$.
If $X'=C_1\cup C_2$ with $C_1\cdot C_2=P$ and we consider on  $X$
(resp. $X'$ a polarization $H$ (resp. $H'$) as above, then there
are isomorphisms
$$\Jac^d(X)_s-\jac^d(X)_s\simeq \Jac^{d-1}(X')_s\quad \text{ and }$$
$$\Jac^d(X)-\jac^d(X)\simeq \Jac^{d-1}(X')\, .$$
\end{proposition}

\begin{proof} Let $\frak{m}$ be the maximal ideal of $2P$ in $X$
and let us denote $\frak{B}=\End_{\mathcal{O}_X}(\frak{m})$ and
$X^*=\Spe(\frak{B})$. Taking into account that $X$ is Gorenstein
and that the natural morphism $\phi\colon X'\to X$ is finite,
birational and non trivial at $2P$, by (2) and (3) of Lemma
\ref{l:gorenstein}, we get a morphism $X'\to X^*$. Since $X'$ and
$X^*$ are both of arithmetic genus zero, this morphism is an
isomorphism. Let $F$ be a pure dimension one sheaf on $X$ of rank
1 and degree $d$ which is not a line bundle, that is, it is not
locally free at the point $2P$. By Lemma \ref{l:nodelinea}, there
is a unique pure dimension one rank 1 sheaf $G$ on $X'$ such that
$\phi_*(G)\simeq F$. Since $\phi$ is finite, birational and
surjective and $\gen(X')=\gen(X)-1$, by Lemma \ref{l:estporimdi},
$G$ has degree $d-1$ and it is (semi)stable with respect to $H'$
if and only if $F$ is (semi)stable with respect to $H$. Therefore,
the direct image $\phi_*$ produces the desired isomorphisms.
\end{proof}

\subsection{The description for the fibers of type $IV$}\label{desX'} In this
part, $X$ will denote a fiber of an elliptic fibration of type
$IV$, that is, $X=C_1\cup C_2\cup C_3$ with $C_i$ rational smooth
curves and $C_1\cdot C_2=C_1\cdot C_3=C_2\cdot C_3=P$ (figure 2).

\hspace{4.6truecm}\includegraphics[scale=0.3]{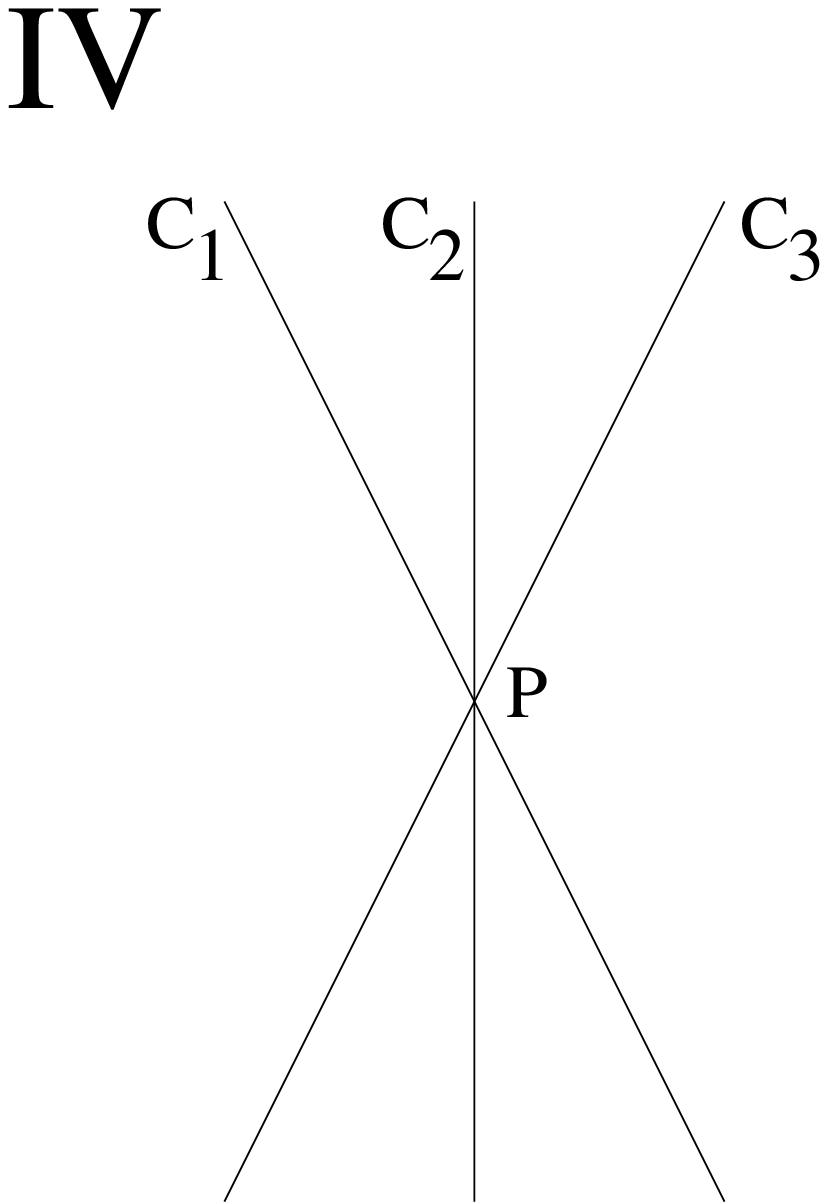}
\vspace{0.2truecm}

\centerline{Figure 2.} \vspace{0.5truecm}

Let $X'=C_1\cup C_2\cup C_3$ be the curve obtained by gluing the
irreducible components of $X$ at the point $P$ in such a way that
this curve $X'$ cannot be embedded, even locally, in a smooth
surface (see figure 3).

\hspace{5truecm}\includegraphics[scale=0.3]{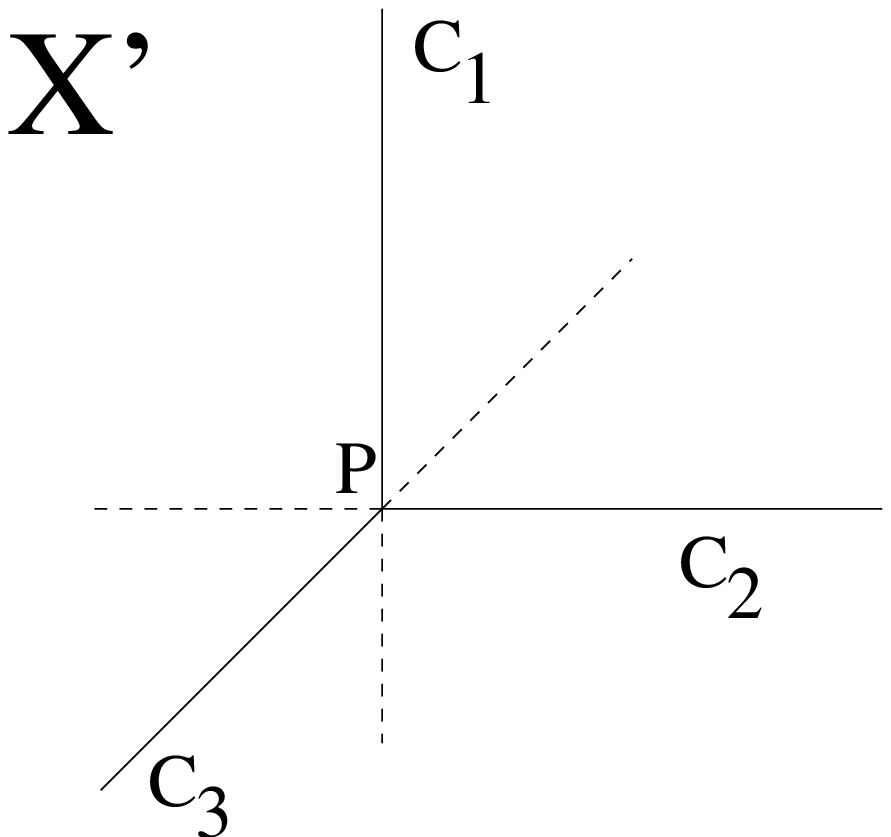}
\vspace{0.2truecm}

\centerline{Figure 3.} \vspace{0.5truecm}

There exists a finite birational morphism $\phi\colon X'\to X$
which is an isomorphism outside $P$ and $\phi^{-1}(P)$ is the
point $P$ of $X'$. Since the completions of the local rings of $X$
and $X'$ at $P$ can be written as
$$\widehat{\mathcal{O}}_{X,P}=\kappa[[x,y]]/xy(y-x)\quad \text{ and
}\quad
\widehat{\mathcal{O}}_{X',P}=\kappa[[\bar{x},\bar{y},\bar{z}]]/(\bar{x}\bar{y},\bar{x}\bar{z},\bar{y}\bar{z})$$
with $x=\bar{x}+\bar{y}$, $y=\bar{y}+\bar{z}$ by $\phi$, it is
easy to prove that $\frak{m}_{X',P}^2\subset\frak{m}_{X,P}$. Then,
from Corollary \ref{c:grupoaditivo}, there is an exact sequence
\begin{equation}\label{e:picIV}
0\to \mathbb{G}_a\to \Pic(X)\xrightarrow{\phi^*}\Pic(X')\to 0\, .
\end{equation}
Moreover, the Picard group of $X'$ is isomorphic to
$\prod_{i=1}^3\Pic(C_i)$. Indeed, if $\tilde{X}=\sqcup_{i=1}^3
C_i$ is the total normalization of $X'$ at $P$ and $\pi\colon
\tilde{X}\to X'$ is the projection map, by Proposition
\ref{p:grothendieck}, we have the following exact sequence:
\begin{equation}\label{e:picX'}
0\to (\pi_*{\mathcal{O}}_{\tilde{X},P}^*/\mathcal{O}_{X',P}^*)/\im
\Coh^0(\tilde{X},\mathcal{O}_{\tilde{X}}^*)\to
\Pic(X')\xrightarrow{\pi^*} \Pic(\tilde{X})\to 0\, .
\end{equation}  Considering the exact sequence
$$0\to \mathcal{O}_{X'}\to \mathcal{O}_{C_1}\oplus \mathcal{O}_{C_2}\oplus
\mathcal{O}_{C_3}\xrightarrow{\beta} \kappa\oplus\kappa\to 0$$
where the morphism $\beta$ is defined by:
$$\beta(s_1,s_2,s_3)=(s_1(P)-s_2(P), s_1(P)-s_3(P))\, ,$$ we get
$$\pi_*{\mathcal{O}}_{\tilde{X},P}^*/\mathcal{O}_{X',P}^*\simeq
\kappa^*\oplus\kappa^*\, .$$ Since the morphism
$\Coh^0(\tilde{X},\mathcal{O}_{\tilde{X}}^*)=\oplus_{i=1}^3\Coh^0(C_i,\mathcal{O}_{C_i}^*)
\xrightarrow{\beta} \kappa^*\oplus\kappa^*$, which is given by
$$\beta(u_1,u_2,u_3)=\big{(}\frac{u_1(P)}{u_2(P)},\frac{u_1(P)}{u_3(P)}\big{)}\,
,$$ is surjective, by \eqref{e:picX'}, we obtain that
$$\Pic(X')\simeq \Pic(\tilde{X})\simeq \prod_{i=1}^3\Pic(C_i)\,
.$$ From \eqref{e:picIV}, we conclude that the Picard group of the
curve $X$ is given by the following exact sequence:
\begin{equation}\label{e:PicIV}
0\to \mathbb{G}_a\to \Pic(X)\to \prod_{i=1}^3\Pic(C_i)\to 0\, .
\end{equation}

We start now with the (semi)stability analysis. Let $H$ be a
polarization on $X$ of degree $h$. As in \eqref{notacion}, if $d$
is the degree of the Jacobian we are considering, $b$ will be the
residue class of $d-1$ modulo $h$. When the numbers
$k_{C_i}=\frac{h_{C_i}(b+1)}{h}$ are not integers, we write
$k_{C_i}=[k_{C_i}]+a_i$ with $0<a_i<1$. Since $\sum_{i=1}^3
k_{C_i}\in \Z$, then $\sum_{i=1}^3 a_i$ is equal to 1 or 2. Thus
the proposition describing the schemes $\jac^d(X)_s$ and
$\jac^d(X)$ is the following:

\begin{proposition}\label{p:lineaIV} Let $X=C_1\cup C_2\cup C_3$ be a curve of type
$IV$. We have the following three cases:

a) If $k_{C_i}\in \Z$ for $i=1,2,3$, then there is an exact
sequence
$$0\to \mathbb{G}_a\to \jac^d(X)_s\to\prod_{i=1}^3\Pic^{h_{C_i}t+k_{C_i}}(C_i)\to 0 \, , \text{ and}$$
$$\jac^d(X)-\jac^d(X)_s=\prod_{i=1}^3\Pic^{d_i}(C_i)$$ where $d_i=h_{C_i}t+k_{C_i}-1$
for all $i$.

b) Since it is not possible to have $k_{C_i}\notin \Z$ only for
one index $i$ because the sum of these three numbers is an
integer, the  following case to considering is $k_{C_i}$ integer
only for one $i$. Set $k_{C_1}\in \Z$ and $k_{C_i}\notin \Z$ for
$i=2,3$. Then, there is an exact sequence $$0\to \mathbb{G}_a\to
\jac^d(X)_s\to K \to 0
$$ where $K=\bigsqcup_{j,l=2}^3\Pic^{h_{C_1}t+k_{C_1}}(C_1)\times
\Pic^{h_{C_j}t+[k_{C_j}]}(C_j)\times
\Pic^{h_{C_l}t+[k_{C_l}]+1}(C_l)$. In this case, we have
$$\jac^d(X)-\jac^d(X)_s=\prod_{i=1}^3\Pic^{d_i}(C_i)$$  where $d_1=h_{C_1}t+k_{C_1}-1$
and $d_2$, $d_3$ are the integer numbers obtained by applying
Algorithm \ref{algoritmo} to the tree-like curve $C_2\cup C_3$ for
a sheaf of degree $d-d_1-2$.

c) Suppose that $k_{C_i}\notin \Z$ for $i=1,2,3$.

1. If $\sum_i a_i=1$, there is an exact sequence
$$0\to \mathbb{G}_a\to\jac^d(X)_s\to K\to 0$$ with $K=\bigsqcup_{i,j,l=1}^3\Pic^{h_{C_i}t+[k_{C_i}]+1}(C_i)\times
\Pic^{h_{C_j}t+[k_{C_j}]}(C_j)\times
\Pic^{h_{C_l}t+[k_{C_l}]}(C_l)$.

2. If $\sum_i a_i=2$, there is an exact sequence
$$0\to \mathbb{G}_a\to \jac^d(X)_s\to K\to 0$$ with
$K=\bigsqcup_{i,j,l=1}^3\Pic^{h_{C_i}t+[k_{C_i}]+1}(C_i)\times
\Pic^{h_{C_j}t+[k_{C_j}]+1}(C_j)\times
\Pic^{h_{C_l}t+[k_{C_l}]}(C_l)$.

In both cases, $\jac^d(X)-\jac^d(X)_s$ is empty, that is, every
semistable line bundle is stable.
\end{proposition}

\begin{proof} Since
$\chi(\mathcal{O}_X)=\chi(\mathcal{O}_{C_1})+\chi(\mathcal{O}_{C_2})+\chi(\mathcal{O}_{C_3})-3$,
for every line bundle $L$ of degree $d$ on $X$, one has an exact
sequence $$0\to L\to L_{C_1}\oplus L_{C_2}\oplus L_{C_3}\to T\to
0$$ where $T$ is a torsion sheaf supported at the point $P$ and
$\chi(T)=3$.

The only connected subcurves of $X$ are $C_i$, $i=1,2,3$, and
their complements $\overline{C_i}=C_j\cup C_l$ and $C_i\cdot
\overline{C_i}=2$ so that, by Lemma
\ref{l:desigualdadesgenerales}, the sheaf $L$ is (semi)stable if
and only if
\begin{equation}\label{d:d2}
h_{C_i}t+k_{C_i}-1\des d_{C_i} \des h_{C_i}t+k_{C_i}+1 ,\text{ for
} i=1,2,3\, .
\end{equation}
Taking into account that $d_{C_i}\in \Z$ and
$d=d_{C_1}+d_{C_2}+d_{C_3}$, in case a) $L$ is stable if and only
if $d_{C_i}=h_{C_i}t+k_{C_i}$ for $i=1,2,3$. Since the following
case is $k_{C_i}\notin \Z$ for two indices, reordering the
irreducible components of $X$, we can assume that we are in case
b). In this case, $L$ is stable if and only if
$d_{C_1}=h_{C_1}t+k_{C_1}$, $d_{C_2}=h_{C_2}t+[k_{C_2}]$ and
$d_{C_3}=h_{C_3}t+[k_{C_3}]+1$ or the same with the roles of $C_2$
and $C_3$ intertwined. Finally, in case c), $L$ is stable if and
only if $d_{C_i}=h_{C_i}t+[k_{C_i}]+\epsilon_i$ where
$\epsilon_i=0\text{ or }1$. Therefore, reordering the irreducible
components of $X$ if it were necessary, the possibilities are:

1. if $\sum_i a_i=1$, then $\epsilon_1=1$ and
$\epsilon_2=\epsilon_3=0$.

2. if $\sum_i a_i=2$, then $\epsilon_1=\epsilon_2=1$ and
$\epsilon_3=0$.

This together with \eqref{e:PicIV} proves the exact sequences of
the statement.

When $L$ is a strictly semistable line bundle on $X$ of degree
$d$, $d_{C_i}$ is equal to one of the two extremal values of
\eqref{d:d2}, and then $k_{C_i}\in \Z$, for some $i=1,2,3$.
Reordering the irreducible components of $X$, we can assume that
$k_{C_1}\in \Z$ and $d_{C_1}=h_{C_1}t+k_{C_1}-1$. Hence,
$\mu_H(L_{C_1})=\mu_H(L^{C_1})=d$ so that, by Lemma
\ref{l:estabilidadinducida}, $L_{C_1}$ and $L^{C_1}$ are
semistable with respect to $H_{C_1}$ and $H_{C_2\cup C_3}$
respectively. Since $C_1$ is an integral curve $L_{C_1}\in
\Pic^{d_{C_1}}(C_1)$ is stable. On the other hand, since $C_2\cup
C_3$ is a tree-like curve, by Theorem \ref{t:treelike}, we have
that $[L^{C_1}]\in \prod_{i=2}^3\Pic^{d_i}(C_i)$ where $d_2$ and
$d_3$ are the integer numbers obtained by applying Algorithm
\ref{algoritmo} to $C_2\cup C_3$ for the sheaf $L^{C_1}$ which has
degree $d-d_{C_1}-2$. Bearing in mind that if $k_{C_i}\in \Z$ for
$i=2,3$, the only final subcurves of $C_2\cup C_3$ are its
irreducible components, we conclude that the $S$-equivalence class
of $L$ belongs to $\prod_{i=1}^3\Pic^{d_i}(C_i)$ with $d_i$ the
integers of the statement.
\end{proof}

We study now the set of boundary points of the connected component
$\Jac^d(X)$.

\begin{proposition}\label{p:bordeIV}
Let $X=C_1\cup C_2\cup C_3$ be a curve of type $IV$ with $C_1\cdot
C_2=C_1\cdot C_3=C_2\cdot C_3=P$ and let $X'=C_1\cup C_2\cup C_3$
be the curve of the figure 3. Let $\phi\colon X'\to X$ denote the
natural morphism. Let $H$ be a polarization on $X$ of degree $h$
such that $H'=\phi^*(H)$ is also of degree $h$. Then, considering
on $X'$ the polarization $H'$, there are isomorphisms
$$\Jac^d(X)_s-\jac^d(X)_s\simeq \Jac^{d-1}(X')_s\quad \text{ and }$$
$$\Jac^d(X)-\jac^d(X)\simeq \Jac^{d-1}(X')\, .$$
\end{proposition}

\begin{proof} Since the curve $X$ is Gorenstein, $\gen(X')=\gen(X)-1$
and $\phi$ is a finite birational morphism non trivial at the
point $P$, arguing as in the proof of Proposition
\ref{p:bordeIII}, the morphism $\phi_*$ produces the desired
isomorphisms.
\end{proof}

Thus the description of the connected component $\Jac^d(X)$ when
$X$ is a fiber of type $IV$ implies the analysis of Simpson
schemes $\jac^d(X')$ and $\Jac^d(X')$ for the curve $X'$ of the
figure 3. This analysis is given in the following subsection which
completes then the description for the curves of type $IV$.

\subsection{The description for the curve $X'$}\label{s:curvaX'}

Let $X'=C_1\cup C_2\cup C_3$ be the curve of the figure 3. Let
$H'$ be a polarization on $X'$ of degree $h$. With the notations
we come using, we have the following

\begin{proposition} If $X'=C_1\cup C_2\cup C_3$ is the curve of
the figure 3, we have the following two cases:

a) If $k_{C_i}\in \Z$ for some $i=1,2,3$, then
$$\jac^d(X')_s=\emptyset \text{ and }$$
$$\jac^d(X')-\jac^d(X')_s=\prod_{i=1}^3\Pic^{d_i}(C_i)$$ where $d_i=h_{C_i}t+k_{C_i}-1$ if $k_{C_i}\in \Z$ and
$d_j$, $d_k$, $j,k\neq i$ are the integer numbers  obtained by
applying Algorithm \ref{algoritmo} to the tree-like curve
$\overline{C_i}$ for a sheaf of degree $d-d_i-1$.

b) Suppose that $k_{C_i}\notin \Z$ for $i=1,2,3$.

1. If $\sum_i a_i=1$, then
$$\jac^d(X')_s=\prod_{i=1}^3\Pic^{h_{C_i}t+[k_{C_i}]}(C_i)\, .$$
2. If $\sum_i a_i=2$, then
$$\jac^d(X')_s=\emptyset\, .$$
In these two cases,
$$\jac^d(X')-\jac^d(X')_s=\emptyset\, .$$
\end{proposition}

\begin{proof} We have seen that the Picard group of this curve
$X'$ is isomorphic to the direct product of the Picard groups of
its irreducible components. Moreover, if $L$ is a line bundle on
$X'$ of degree $d$, by Lemma \ref{l:desigualdadesgenerales}, $L$
is (semi)stable if and only if $$h_{C_i}t+k_{C_i}-1\des d_{C_i}
\des h_{C_i}t+k_{C_i} ,\text{ for } i=1,2,3\, .$$ Then, the result
follows arguing as in the proof of Proposition  \ref{p:lineaIV}.
\end{proof}

In order to determine boundary points of $\Jac^d(X')$, we need the
following lemma that describes the stalk of a pure dimension one
sheaf on the curve $X'$ at its only singular point $P$.

\begin{lemma}\label{l:localizacion2}
If $F$ is a pure dimension one sheaf on the curve $X'$, then
$$F_{P}\simeq \mathcal{O}_{X',P}^{a}\oplus \mathcal{O}_{C_{12},P}^{a_{12}} \oplus \mathcal{O}_{C_{13},P}^{a_{13}}\oplus
\mathcal{O}_{C_{23},P}^{a_{23}}\oplus
\mathcal{O}_{C_1,P}^{a_1}\oplus \mathcal{O}_{C_2,P}^{a_2}\oplus
\mathcal{O}_{C_3,P}^{a_3}$$ where $C_{ij}=C_i\cup C_j$ and $a$,
$a_{ij}$ and $a_i$ are integer numbers determined by the following
equalities:
\begin{align}
&a+a_{12}+a_{13}+a_1=\rg(F_{P}\underset{\mathcal{O}_{X',
P}}\otimes
\mathcal{O}_{C_1,P})\notag\\
&a+a_{12}+a_{23}+a_2=\rg(F_{P}\underset{\mathcal{O}_{X',
P}}\otimes
\mathcal{O}_{C_2,P})\notag\\
&a+a_{13}+a_{23}+a_3=\rg(F_{P}\underset{\mathcal{O}_{X',
P}}\otimes \mathcal{O}_{C_3,P})\notag\\
&a+a_{12}+\frac{a_{13}}{2}+\frac{a_{23}}{2}+\frac{a_1}{2}+\frac{a_2}{2}=\rg(F_{P}\underset{\mathcal{O}_{X',
P}}\otimes \mathcal{O}_{C_{12},P})\notag\\
&a+\frac{a_{12}}{2}+a_{13}+\frac{a_{23}}{2}+\frac{a_1}{2}+\frac{a_3}{2}=\rg(F_{P}\underset{\mathcal{O}_{X',
P}}\otimes \mathcal{O}_{C_{13},P})\notag\\
&a+\frac{a_{12}}{2}+\frac{a_{13}}{2}+a_{23}+\frac{a_2}{2}+\frac{a_3}{2}=\rg(F_{P}\underset{\mathcal{O}_{X',
P}}\otimes \mathcal{O}_{C_{23},P})\notag\\
&a+a_{12}+a_{13}+a_{23}+a_1+a_2+a_3=\rg(F_{P}\otimes \kappa)\, .
\notag
\end{align}
\end{lemma}

\begin{proof} Since $F$ is a pure dimension one sheaf, the $\mathcal{O}_{X',P}$-module
$M:=F_P$ has depth 1. For $i=1,2,3$, let $t_i$ be a local
parameter of $C_i$ at $P$. If $\frak{m}$ denotes the maximal ideal
of $P$ in  $X'$ and $\frak{m}_i$ is the maximal ideal of $P$ in
$C_i$, we have $\frak{m}=\oplus_{i=1}^3 \frak{m}_i$ so that
$t_it_j=0$ for $i\neq j$. Since $t_iM$ is a torsion free
$\mathcal{O}_{C_i,P}$-module, it is free:
\begin{equation}\label{e:isomorfismos}
t_iM\simeq \mathcal{O}_{C_i,P}^{r_i}, \quad i=1,2,3.
\end{equation}
The map
\begin{align*}
 \Psi\colon M&\to t_1M\oplus t_2M\oplus t_3M \\ m&\mapsto
 (t_1m,t_2m,t_3m)
\end{align*} is injective because $M$ has depth 1.
Consider the inclusion map
\begin{align*}
i\colon t_1M\oplus t_2M\oplus t_3M &\hookrightarrow M\\
(t_1m,t_2m',t_3m'')&\mapsto t_1m+t_2m'+t_3m''\, .
\end{align*}
From \eqref{e:isomorfismos}, one has that
$$\Psi(\im(i))=\frak{m}_1^{r_1}\oplus
\frak{m}_2^{r_2}\oplus \frak{m}_3^{r_3}$$ and then
$$(t_1M\oplus t_2M\oplus t_3M)/\Psi(\im(i))=\kappa^{r_1}\oplus
\kappa^{r_2}\oplus \kappa^{r_3}\, .$$ Hence, there is a map
$$\chi\colon M\to \kappa^{r_1}\oplus
\kappa^{r_2}\oplus \kappa^{r_3}$$ whose kernel is $t_1M\oplus
t_2M\oplus t_3M$.

For $i=1,2,3$, let $M_i$ be the following submodule of $M$:
$$M_i=\{ f\in M\text{ such that } f\notin t_iM \text{ and }
t_jf=t_kf=0 \text{ for } j,k\neq i\}\, .$$ If
$N_i:=\chi(M_i)\subseteq \kappa^{r_i}$, then $N_i=\kappa^{r_i}\cap
\im(\chi)$. Indeed, let $u$ be a nonzero element of
$\kappa^{r_i}\cap \im(\chi)$ with $u=\chi(f)$ for $f\in M$. Since
$t_jf\in t_j^2M$ for $j\neq i$, we can write $t_jf=t_j^2g$ for
some $g\in M$. If $f'=f-t_jg$, then $\chi(f)=\chi(f')=u$ and
$f'\in M_i$ which proves the claim.

Let $f_1,\hdots,f_n$ be elements of $M_1$ such that
$\chi(f_1),\hdots,\chi(f_n)$ are free and let us check that
$f_1,\hdots,f_n$ are free over $\mathcal{O}_{C_1,P}$. Suppose that
$$\alpha_1f_1+\hdots+\alpha_nf_n=0$$
with $\alpha_i\in \mathcal{O}_{C_1,P}$ and $\alpha_i\neq 0$ for
some $i$. Let $t_1^m$ be the maximal power of $t_1$ dividing all
$\alpha_i$ and write $\alpha_i=t_1^m\beta_i$ for $i=1,\hdots n$.
Since $t_it_j=0$ for $i\neq j$, we have that
$$t_i(t_1^{m-1}\beta_1f_1+\hdots +t_1^{m-1}\beta_nf_n)=0$$ for $i=1,2,3$,
which implies that
$$t_1^{m-1}\beta_1f_1+\hdots +t_1^{m-1}\beta_nf_n=0$$ because $M$
has depth 1.

Recurrently, we get
$$\beta_1f_1+\hdots +\beta_nf_n=0\, .$$ Therefore,
$$\beta_1(P)\chi(f_1)+\hdots +\beta(P)\chi(f_n)=0$$ which is
absurd because $\beta_i(P)\neq 0$ for some $i=1,\hdots, n$ and
$\chi(f_1)\hdots, \chi(f_n)$ are free. Thus $M_1$ is a free
$\mathcal{O}_{C_1,P}$-module.

The same argument proves that $M_2$ (resp. $M_3$) is a free
$\mathcal{O}_{C_2,P}$ (resp. $\mathcal{O}_{C_3,P}$) module.

One has $M_i\cap M_j=\{ 0\}$ for $i\neq j$ as $t_k(M_i\cap M_j)=0$
for $k=1,2,3$ and $M$ has depth 1.

Let $K_1$ be a vector subspace of $\im(\chi)$ which is
supplementary of $N_1\oplus N_2\oplus N_3$. For $i\neq j$ with
$i,j=1,2,3$, let $M_{ij}$ the following submodule of $M$:
$$M_{ij}=\{ f\in M\text{ such that } f\notin t_iM\oplus t_jM , \ t_kf=0 \text{ for }k\neq i,j, \chi(f)\in
K_1\}\, .$$ Arguing as above, if $N_{ij}:=\chi(M_{ij})$, then
$N_{ij}=(\kappa^{r_i}\oplus \kappa^{r_j})\cap \im(\chi)$.

Let $f_1\hdots,f_n$ be elements of $M_{12}$ such that
$\chi(f_1),\hdots,\chi(f_n)$ are free and let us prove that
$f_1,\hdots,f_n$ are also free over  $\mathcal{O}_{C_{12},P}$.
Suppose that
\begin{equation}\label{e:igualdad1}
\alpha_1f_1+\hdots +\alpha_nf_n=0 \end{equation} where
$\alpha_i\in \mathcal{O}_{C_{12},P}$ are not all nulls. Since
$\chi(f_1),\hdots, \chi(f_n)$ are free, $\alpha_i(P)=0$ for
$i=1\hdots,n$. Let $t_1^m$ and $t_2^s$ be the maximal powers of
$t_1$ and $t_2$ dividing all $\alpha_i$ and write $\alpha_i=t_1^m
u_i+t_2^s v_i$, $i=1,\hdots,n$. Since $t_1t_2=0$, we have that
either $u_i(P)$ are not all nulls or  $v_i(P)$ are not all nulls.
Suppose that $u_i(P)$ are not all nulls. By multiplying by $t_1$
the equality \eqref{e:igualdad1}, we get
$$t_1(t_1^mu_1f_1+\hdots+t_1^mu_nf_n)=0\, .$$ Since $t_i(t_1^mu_1f_1+\hdots+t_1^mu_nf_n)=0$ for $i=2,3$ as well
and $M$ has depth 1, we have that
$$t_1^{m}u_1f_1+\hdots t_1^mu_nf_n=0\, .$$ Recurrently, we get
$$t_1(u_1f_1+\hdots+u_nf_n)=0\, .$$
Set $w:=u_1f_1+\hdots+u_nf_n$. Since $\omega\notin M_2$,  then
$\omega\in t_2M$ so that $\chi(\omega)$ is zero, but this is
absurd because $u_i(P)$ are not all nulls and  $\chi(f_1),\hdots,
\chi(f_n)$ are free. Thus the $\mathcal{O}_{C_{12}}$-module
$M_{12}$ is free.

Analogously, $M_{13}$ and $M_{23}$ are free modules and it is not
difficult to prove that all intersections of $M_{ij}$ and $M_i$
are zero.

Let $K_2$ a vector subspace of $\im(\chi)$ supplementary of
$N_1\oplus N_2\oplus N_3\oplus N_{12}\oplus N_{13}\oplus  N_{23}$.

Let $g_1,\hdots,g_d$ be elements of $M$ such that
$\chi(g_1),\hdots,\chi(g_d)$ are a basis of $K_2$. Let $M_0$ be
the submodule of $M$ generated by $g_1,\hdots, g_d$. Let us prove
that $g_1,\hdots,g_d$ are free over $\mathcal{O}_{X',P}$. Suppose
that
\begin{equation}\label{e:igualdad2}
\alpha_1g_1+\hdots +\alpha_df_d=0 \end{equation} where
$\alpha_i\in \mathcal{O}_{X',P}$ are not all nulls. Since
$\chi(g_1),\hdots,\chi(g_d)$ are free, $\alpha_i(P)=0$ for
$i=1,\hdots,d$.  Let $t_1^m$, $t_2^s$ and $t_3^r$ be the maximal
powers of $t_1$, $t_2$ and $t_3$ dividing  $\alpha_i$ for
$i=1,\hdots,d$ and write $\alpha_i=t_1^m u_i+t_2^s v_i+t_3^rw_i$,
$i=1,\hdots,d$. Since $t_it_j=0$ for $i\neq j$, one of the
following conditions holds: $u_i(P)\neq 0$ for some $i$,
$v_i(P)\neq 0$ for some $i$ or $w_i(P)\neq 0$ for some $i$.
Suppose that $u_i(P)\neq 0$ for some $i$. By multiplying by $t_1$
equality \eqref{e:igualdad2} and arguing as in the former case, we
get
$$t_1(u_1g_1+\hdots u_dg_d)=0\, .$$ If $w':=u_1g_1+\hdots+u_dg_d$, then
$\chi(w')$ is a nonzero element of $K_2$. But this is impossible
because $w'\in t_2M\oplus t_3M$ and  $t_iw'\neq 0$ for $i=1,2,3$.
Hence $M_0$ is a free $\mathcal{O}_{X',P}$-module.

One has also $M_0\cap M_i=M_0\cap M_{ij}=\{ 0\}$ for all
$i,j=1,2,3$.

From Nakayama's Lemma, we conclude then that $$M=M_0\oplus
M_1\oplus M_2\oplus M_3\oplus M_{12}\oplus M_{13}\oplus M_{23}$$
which proves the result.
\end{proof}

\begin{proposition}\label{p:bordex'} If $X'=C_1\cup C_2\cup C_3$ is the curve of the figure
3, then the only stable pure dimension one rank 1 sheaves on $X'$
are stable line bundles, that is,
$$\Jac^d(X')_s-\jac^d(X')_s=\emptyset\, .$$

Moreover, we have that

a) if  $k_{C_i}\notin \Z$ for all $i=1,2,3$, then the set of
boundary points is empty, that is,
$$\Jac^d(X')-\jac^d(X')=\emptyset\, ,$$

b) if $k_{C_i}\in \Z$ for some $i=1,2,3$, then
$$\Jac^d(X')-\jac^d(X')=\prod_{i=1}^3\Pic^{d_i}(C_i)$$ where
$d_i=h_{C_i}t+k_{C_i}-1$ and $d_j$ $j\neq i$ are the integer
numbers obtained by applying Algorithm \ref{algoritmo} to the
tree-like curve $\overline{C_i}$ for a sheaf of degree $d-d_i-1$.
\end{proposition}

\begin{proof} Let $F$ be a pure dimension one sheaf on $X'$ of
rank 1 and degree $d$ with respect to the polarization $H'$. If
$F$ is not a line bundle, from Lemma \ref{l:localizacion2}, we get
that, reordering the irreducible components of $X'$ if it were
necessary, the sheaf $F$ is isomorphic either to $F_{C_1}\oplus
F_{C_2}\oplus F_{C_3}$ or to $F_{C_1\cup C_2}\oplus F_{C_3}$. So
$F$ is not stable.

Suppose that $F$ is semistable and let us determine  its
$S$-equivalence class. If $F\simeq F_{C_1}\oplus F_{C_2}\oplus
F_{C_3}$, the condition $\mu_{H}(F_{C_i})=d$ for $i=1,2,3$ implies
that $d_{C_i}=h_{C_i}t+k_{C_i}-1$. Then $k_{C_i}\in \Z$ for
$i=1,2,3$ and the $S$-equivalence class of  $F$ belongs to
$\prod_{i=1}^3 \Pic^{d_i}(C_i)$ where $d_i$ are the integers of
the statement because in this case the only final subcurves are
the irreducible components. If $F\simeq F_{C_1\cup C_2}\oplus
F_{C_3}$, since $\mu_{H}(F_{C_1\cup C_2})=\mu_{H}(F_{C_3})=d$,
then $d_{C_3}=h_{C_3}t+k_{C_3}-1$ so that $k_{C_3}\in \Z$ and the
sheaf $F_{C_1\cup C_2}$ is semistable of degree $d-d_3-1$ with
respect to $H_{C_1\cup C_2}$ on the tree-like curve $C_1\cup C_2$.
By Theorem \ref{t:treelike}, we conclude that $[F]\in
\prod_{i=1}^3\Pic^{d_i}(C_i)$ with $d_i$ the integers of the
statement as well.
\end{proof}

\subsection{The description for the fibers of type $I_N$}
In all this subsection $X$ will be a fiber of an elliptic
fibration of type $I_N$ (figure 4), that is, if $N>2$, then
$X=C_1\cup C_2\cup\hdots \cup C_N$ with $C_1\cdot C_2=C_2\cdot
C_3=\hdots =C_{N-1}\cdot C_N=C_N\cdot C_1=1$, and if $N=2$, then
$X=C_1\cup C_2$ with $C_1\cdot C_2=P+Q$. In both cases, the
irreducible components of $X$ are rational smooth curves.

\includegraphics[scale=0.3]{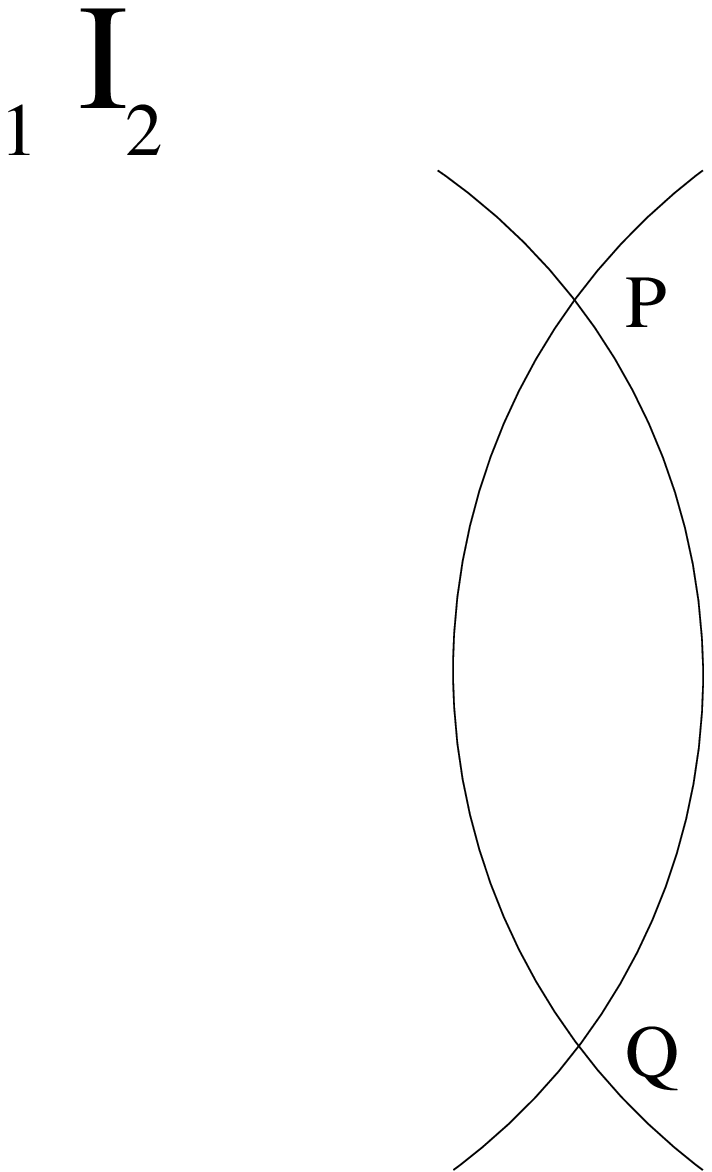}
\hspace{5.5truecm} \includegraphics[scale=0.3]{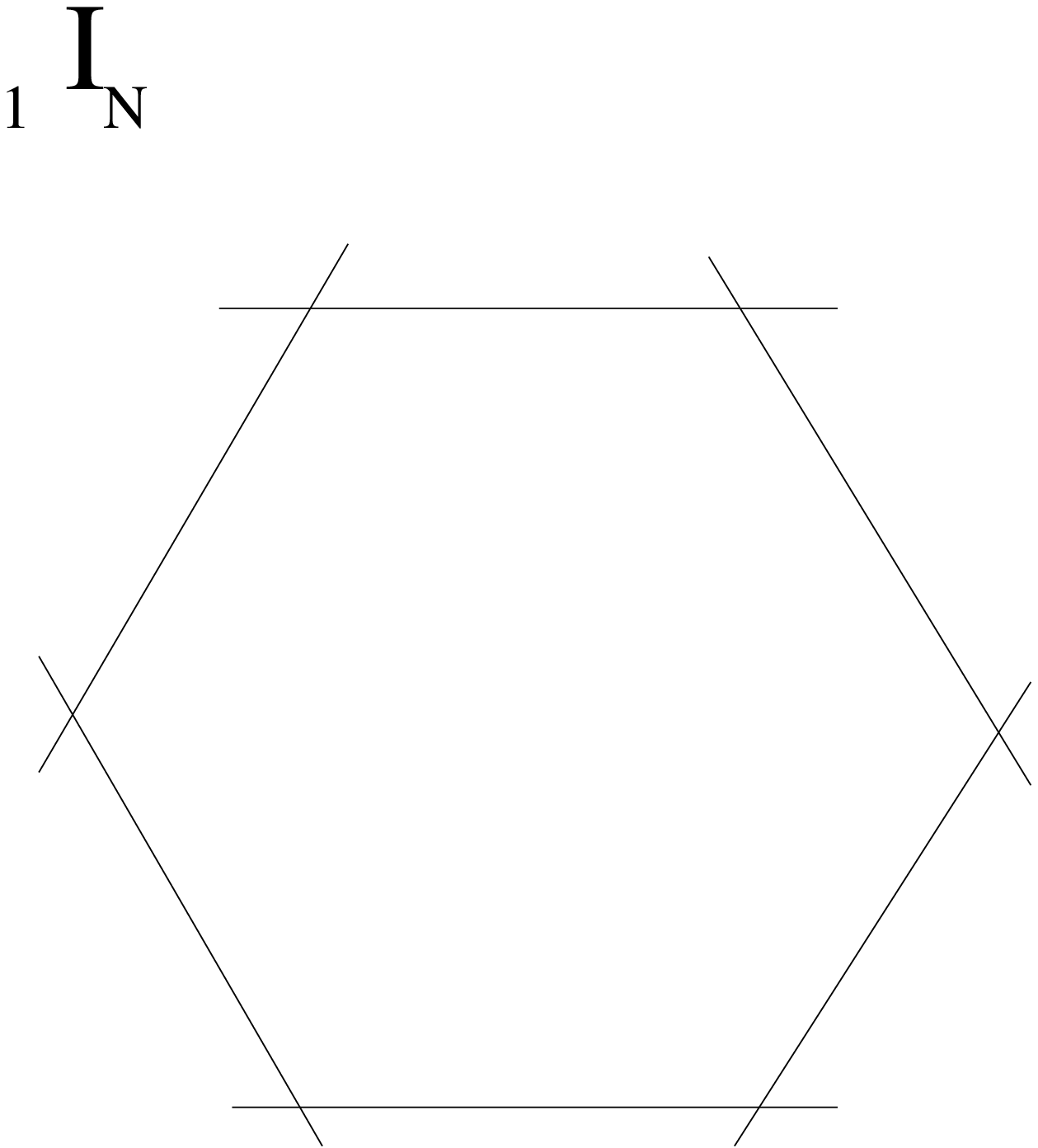}

\centerline{Figure 4.} \vspace{0.5truecm}

Since the number of intersection points of the irreducible
components of $X$ is equal to the number of irreducible
components, by Corollary \ref{c:odaseshadri}, there is an exact
sequence
\begin{equation}\label{e:picIN}0\to
\kappa^*\to \Pic(X)\to \prod_{i=1}^N\Pic(C_i)\to 0\, .
\end{equation}

Let $H$ be a polarization on $X$ of degree $h$. Fixing the degree
$d$ and using the above notations, that is, $b$ is the residue
class of $d-1$ module $h$ and $k_{C_i}=\frac{h_{C_i}(b+1)}{h}$,
let us write $k_{C_i}=[k_{C_i}]+a_i$ with $0<a_i<1$ when these
numbers are not integers and set us $\epsilon_i=0\text{ or }1$.
Then the proposition describing the structure of the Simpson
Jacobian of $X$ of degree $d$ is the following:

\begin{proposition}\label{p:lineaestableIN}
Let $X=C_1\cup \hdots \cup C_N$ be a curve of type $I_N$, $N\geq
2$.

a) If $k_{C_i}\in \Z$ for $i=1,\hdots, N$, there is an exact
sequence
$$0\to \kappa^*\to\jac^d(X)_s\to \prod_{i=1}^{N}\Pic^{h_{C_i}t+k_{C_i}}(C_i)\to 0\, .$$

b) If there are exactly  $r\geq 2$ indices, say $i_1,\hdots, i_r$,
such that the numbers $k_{C_{i_1}},\hdots ,k_{C_{i_r}}$ are not
integers, then $\jac^d(X)_s$ is not empty if and only if for any
subset $J\subseteq\{ i_1,\hdots,i_r\}$ such that either
$\cup_{j\in J}C_j$ is a connected subcurve or it is connected by
adding irreducible components $C_i$ with $i\neq i_1, \hdots,i_r$
the following inequalities hold:
$$\sum_{j\in J}\epsilon_{i_j}-1<\sum_{j\in J}a_{i_j}<\sum_{j\in
J}\epsilon_{i_j}+1\, .$$  In this case, there is an exact sequence
$$0\to \kappa^*\to\jac^d(X)_s\to \prod_{s\neq i_j}\Pic^{h_{C_s}t+k_{C_s}}(C_s)\times \prod_{j=1}^{r}
\Pic^{h_{C_{i_j}}t+[k_{C_{i_j}}]+\epsilon_{i_j}}(C_{i_j})\to 0\,
.$$
\end{proposition}

\begin{proof} Since
$\chi(\mathcal{O}_X)=\sum_{i=1}^N\chi(\mathcal{O}_{C_i})-N$, for
every line bundle $L$ on $X$ of degree $d$ there is an exact
sequence $$0\to L\to L_{C_1}\oplus\hdots\oplus L_{C_N}\to T\to 0$$
where $T$ is a torsion sheaf with $\chi(T)=N$. By the exact
sequence \eqref{e:picIN}, it is enough then to find the values of
the degrees $d_{C_i}$ so that $L$ is stable. Note that every
connected subcurve $D$ of $X$ is a tree-like curve of arithmetic
genus 0 and $D\cdot \overline{D}=2$. Hence, by Lemma
\ref{l:desigualdadesgenerales}, $L$ is (semi)stable if and only if
\begin{equation}\label{desIN}
-h_Dt+k_D-1\des d_D\des h_Dt+k_D+1
\end{equation} for every connected subcurve $D$ of $X$.  Since the
degrees $d_{C_i}$ are integers, we have that

a) if $k_{C_i}\in \mathbb{Z}$ for $i=1,\hdots,N$, then $L$ is
stable if and only if  $d_{C_i}=h_{C_i}t+k_{C_i}$ for
$i=1,\hdots,N$, because  $d_D=\sum_j d_{C_j}$ if $D=\cup_j C_j$.

b) Suppose that only $k_{C_{i_1}},\hdots ,k_{C_{i_r}}$, $r\geq 2$,
are not integers. If $L$ is stable, by \eqref{desIN},
$d_{C_s}=h_{C_s}t+k_{C_s}$ for $s\neq i_1,\hdots, i_r$ and
$d_{C_{i_j}}=h_{C_{i_j}}t+[k_{C_{i_j}}]+\epsilon_{i_j}$ with
$\epsilon_{i_j}=0 \text{ or } 1$ for $j=1,\hdots,r$. Conversely,
suppose that $L$ is a line bundle on $X$ obtained by gluing line
bundles $L_i$ on $C_i$ of degrees $d_{C_s}=h_{C_s}t+k_{C_s}$ for
$s\neq i_1,\hdots, i_r$ and
$d_{C_{i_j}}=h_{C_{i_j}}t+[k_{C_{i_j}}]+\epsilon_{i_j}$ with
$\epsilon_{i_j}=0 \text{ or } 1$ for $j=1,\hdots,r$. Let $D$ be a
connected subcurve of $X$. If $D$ contains no component $C_{i_j}$,
one easily check that $d_D$ holds \eqref{desIN}. Otherwise,
consider $J=\{ i_j, j=1,\hdots, r \text{ such that }
C_{i_j}\subset D\}$. Since $J$ satisfies the second condition of
the statement, by the hypothesis, $d_D=h_Dt+k_D-\sum_{j\in
J}a_{i_j}+\sum_{j\in J}\epsilon_{i_j}$ holds \eqref{desIN} as
well. Then $L$ is stable and the proof is complete.
\end{proof}

For the subscheme of strictly semistable line bundles, we have:

\begin{proposition} \label{p:lineasemiIN}
Let $X=C_1\cup \hdots \cup C_N$ be a curve of type $I_N$, $N\geq
2$. Then,
$$\jac^d(X)-\jac^d(X)_s=\prod_{i}\Pic^{d_i}(C_i)\times
\prod_j\Pic^{d_j}(C_j)$$ where $i$ (resp. $j$) runs through the
irreducible components of a connected subcurve $D\subset X$ (resp.
$\overline{D}$) such that $k_D\in \mathbb{Z}$ and $d_i$ (resp.
$d_j$) are the integers obtained by applying Algorithm
\ref{algoritmo} to $D$ (resp. $\overline{D}$) for a sheaf of
degree $h_Dt+k_D-1$ (resp.
$h_{\overline{D}}t+k_{\overline{D}}-1$). If there is not so
connected subcurve $D$ of $X$, then
$\jac^d(X)-\jac^d(X)_s=\emptyset$.
\end{proposition}

\begin{proof} If $L$ is a strictly semistable line bundle on $X$
of degree $d$, there is a connected subcurve $D$ of $X$ such that
$d_D$ is equal to one of the two extremal values of \eqref{desIN}.
In particular, $k_{D}\in \Z$. Suppose that $d_D=h_Dt+k_D-1$ (the
other case is similar). Then $\mu_H(L^D)=\mu_H(L_D)=\mu_H(L)$ and,
by Lemma \ref{l:estabilidadinducida}, the sheaves $L^D$ and $L_D$
are semistable with respect to $H_{\overline D}$ and $H_D$
respectively. Since $D$ and $\overline{D}$ are tree-like curves,
we conclude thanks to Theorem \ref{t:treelike}.
\end{proof}

Arguing as in the proof of Proposition \ref{p:bordeIII}, we get
the following proposition which together with Theorem
\ref{t:treelike} gives us the structure of the border of
$\Jac^d(X)$ when $X$ is a curve of type $I_N$.

\begin{proposition}\label{p:bordeIN}
Let $X=C_1\cup \hdots \cup C_N$ be a curve of type $I_N$, $N\geq
2$ and let $X'=C_1\cup \hdots \cup C_N$ be the tree-like curve
obtained by blowing up $X$ at one of its singular points. Let
$\phi\colon X'\to X$ denote the natural projection. Let $H$ be a
polarization on $X$ of degree $h$ such that $H':=\phi^*H$ is also
of degree $h$. Then, considering on $X'$ the polarization $H'$,
there are isomorphisms
$$\Jac^d(X)_s-\jac^d(X)_s\simeq\Jac^{d-1}(X')_s \, ,$$
$$\Jac^d(X)-\jac^d(X)\simeq\Jac^{d-1}(X')\, .$$
\end{proposition}

It is known that  the pull-black of a stable torsion free sheaf by
a finite morphism of integral curves is stable. Using the above
descriptions, the following example shows that this result is not
longer true for reducible curves.

\begin{example}\label{e:imageninversa}{\rm Let $X=C_1\cup C_2$ be
a curve of type $I_2$ with $C_1\cdot C_2=P+Q$. Let $X'=C_1\cup
C_2$ be the blow-up of $X$ at one of its singular points, say $P$,
and let $\phi\colon X'\to X$ be the natural morphism.

Let $H$ be a polarization on $X$ of degree $h$ such that
$H'=\phi^*(H)$ is also of degree $h$. Let $L_1$ (resp. $L_2$) be a
line bundle on  $C_1$ (resp. $C_2$) of degree $d_1=h_{C_1}(t+1)$
(resp. $d_2=h_{C_2}(t+1)-1$) for some $t \in \Z$. Let $L$ be a
line bundle on $X$ such that $L_{C_i}\simeq L_i$ for $i=1,2$. With
the above notations, since $L$ has degree $d=ht+(h-1)$, we have
that the residue class of $d-\gen(X)$ module $h$ is equal to
$h-2$, $k_{C_i}=h_{C_i}-\frac{h_{C_i}}{h}\notin \Z$ for $i=1,2$
and $a_1=\frac {h_{C_2}}{h}$, $a_2=\frac{h_{C_1}}{h}$,
$\epsilon_1=1$ and $\epsilon_2=0$. Since the conditions of the
Proposition \ref{p:lineaestableIN} hold for all $J\subseteq \{
1,2\}$, $L$ is stable with respect to $H$. The line bundle
$\phi^*(L)$ on $X'$ has degree $d$ with respect to $H'$. The curve
$X'$ is a tree-like curve of arithmetic genus zero, its
irreducible components are ordered according to Lemma
\ref{l:orden} and $X'_1=C_1$. Since now the residue class of
$d-\gen(X')$ is equal to $h-1$ and the corresponding number
$k'_{X'_1}=h_{C_1}$ is integer, by Theorem \ref{t:treelike}, the
sheaf $\phi^*(L)$ is not stable.}
\end{example}

\section{The case of degree zero.}
In this section we want to make a special emphasis in the case of
degree zero. This case is particularly significative because, as
Corollary \ref{c:nopolarizacion} proves, if $X$ is a fiber of type
$III$, $IV$ or $I_N$ of an elliptic fibration, for degree zero the
conditions of semistability are independent of the polarization on
$X$. Moreover, the results of the previous section now take a
simpler and more explicit form.

We need previously the following

\begin{lemma}\label{l:noest} If  $X=C_1\cup \hdots \cup C_N$ is a
polarized tree-like curve whose irreducible components are
rational, then there is no stable pure dimension one sheaf on $X$
of rank 1 and degree -1.
\end{lemma}

\begin{proof} Suppose that the irreducible components of $X$ are
ordered according to Lemma \ref{l:orden}. Since the arithmetic
genus of $X$ is 0 and $d=-1$, using the notations of
\eqref{notacion}, we have that $b=h-1$. Then, $k_{X_i}=h_{X_i}\in
\Z$ for all $i$ and the result follows from Theorem
\ref{t:treelike}.
\end{proof}

\begin{proposition} \label{p:grado0} Let $X=\cup_i C_i$ be a polarized curve of type $III$, $IV$
or $I_N$. Then the following statements hold:
\begin{enumerate}
\item If $L$ is a line bundle on $X$ of degree 0, then $L$ is stable if
and only if $L_{C_i}=\mathcal{O}_{\mathbb{P}^1}$ for all $i$.
\item Every stable pure dimension one sheaf on $X$ of rank 1 and degree 0 is a line bundle.
\end{enumerate}
\end{proposition}

\begin{proof} 1. Since $d=0$ and $\gen(X)=1$, then $b=h-1$, $t=-1$
and $k_{C_i}=h_{C_i}$ is integer for all $i$. Thus the result
follows from Propositions \ref{l:lineaIII}, \ref{p:lineaIV} and
\ref{p:lineaestableIN}.

2. If $F$ is a stable pure dimension one sheaf on $X$ of rank 1
and degree 0 which is not a line bundle, by Propositions
\ref{p:bordeIII}, \ref{p:bordeIV} and \ref{p:bordeIN},
$F=\phi_*(G)$ where $G$ is a stable pure dimension one sheaf of
rank 1 and degree -1 on a curve $X'$ we have determined. When $X$
is of type $III$ or $I_N$, $X'$ is a tree-like curve with rational
components. Then this is impossible by Lemma \ref{l:noest}. When
$X$ is of type $IV$, $X'$ is the curve of the figure 3 and this is
also absurd by the description given in \ref{desX'}.
\end{proof}

\begin{proposition}\label{p:propiasemi} Let $X$ be a polarized curve of type $III$, $IV$ or $I_N$ and
let $L$ be a line bundle on $X$ of degree 0.
\begin{enumerate}
\item If $X$ is of type $III$, $L$ is strictly semistable if and only if
$L_{C_1}=\mathcal{O}_{\mathbb{P}^1}(-1)$ and
$L_{C_2}=\mathcal{O}_{\mathbb{P}^1}(1)$.
\item If $X$ is of type  $IV$, $L$ is strictly semistable if and only if
 $L_{C_1}=\mathcal{O}_{\mathbb{P}^1}(-1)$, $L_{C_2}=\mathcal{O}_{\mathbb{P}^1}$ and $L_{C_3}=\mathcal{O}_{\mathbb{P}^1}(1)$.
\item If $X$ is of type $I_N$, $L$  is strictly semistable if and only if
$L_{C_i}=\mathcal{O}_{\mathbb{P}^1}(r)$ where $r=-1,0$ or $1$ in
such a way that when we remove the components $C_i$ for which
$r=0$ there are neither two consecutive $r=1$ nor two consecutive
$r=-1$.
\end{enumerate}
\end{proposition}

\begin{proof} If $h$ is the degree of the polarization on $X$,
then the residue class of $\dg_H(L)-\gen(X)=-1$ module $h$ is
$b=h-1$ and $t=-1$. Then, for every connected subcurve  $D$ of
$X$, it is $k_{D}=\frac{h_{D}(b+1)}{h}=h_D\in \Z$. Since the
arithmetic genus of $D$ is 0 and $D\cdot \overline{D}=2$, by Lemma
\ref{l:desigualdadesgenerales}, $L$ is strictly (semi)stable if
and only if for every connected $D\subset X$ we have
$$-1\des d_D\des 1$$ and $d_D$ is equal to one of the two extremal
values for some $D$. The result is now straightforward.
\end{proof}

\begin{example}{\rm The only possibilities for a strictly semistable line bundle of degree 0 on a curve of type $I_3$
and $I_4$ are:

\

\includegraphics[scale=0.3]{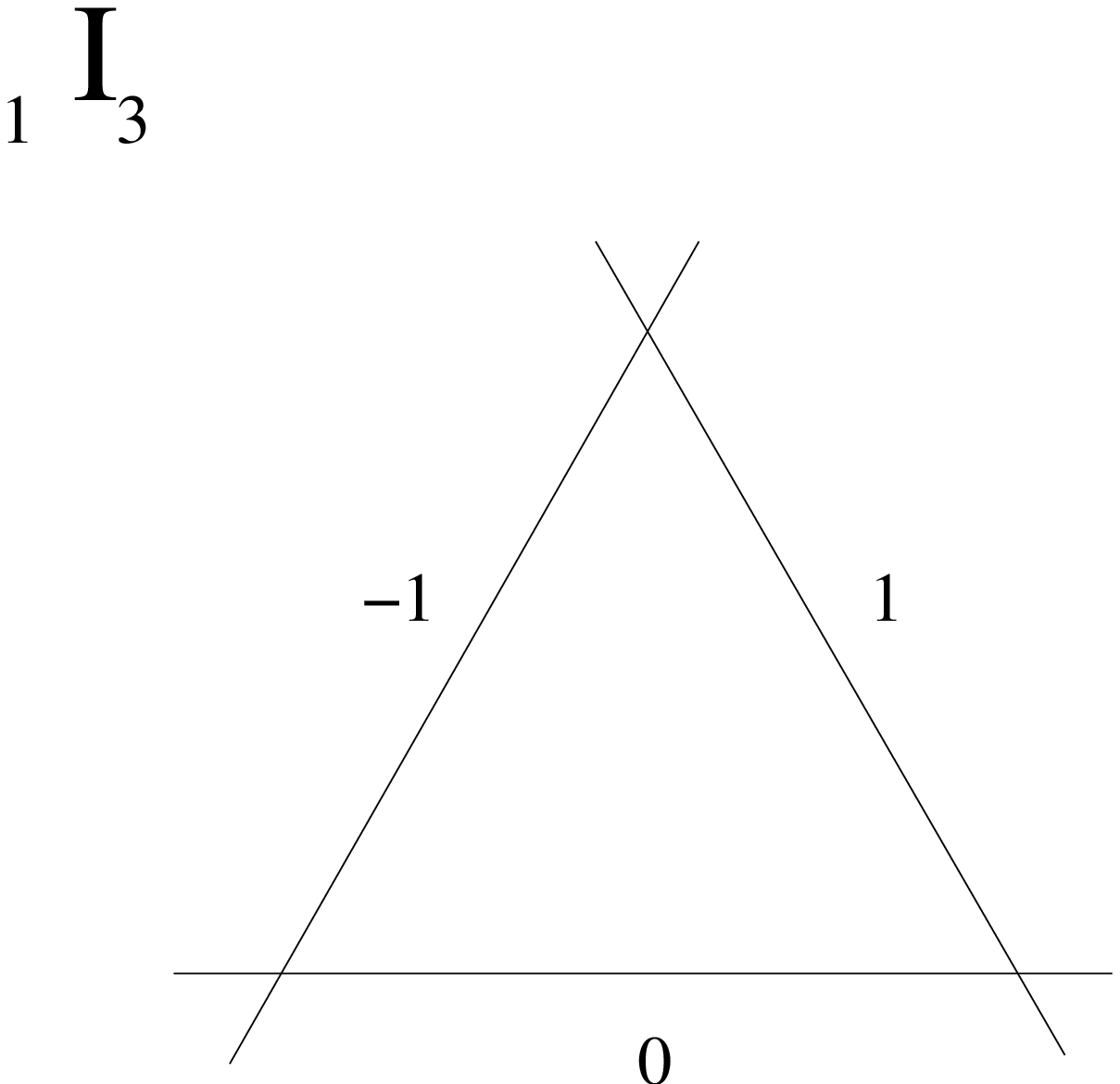}
\hspace{3truecm} \includegraphics[scale=0.3]{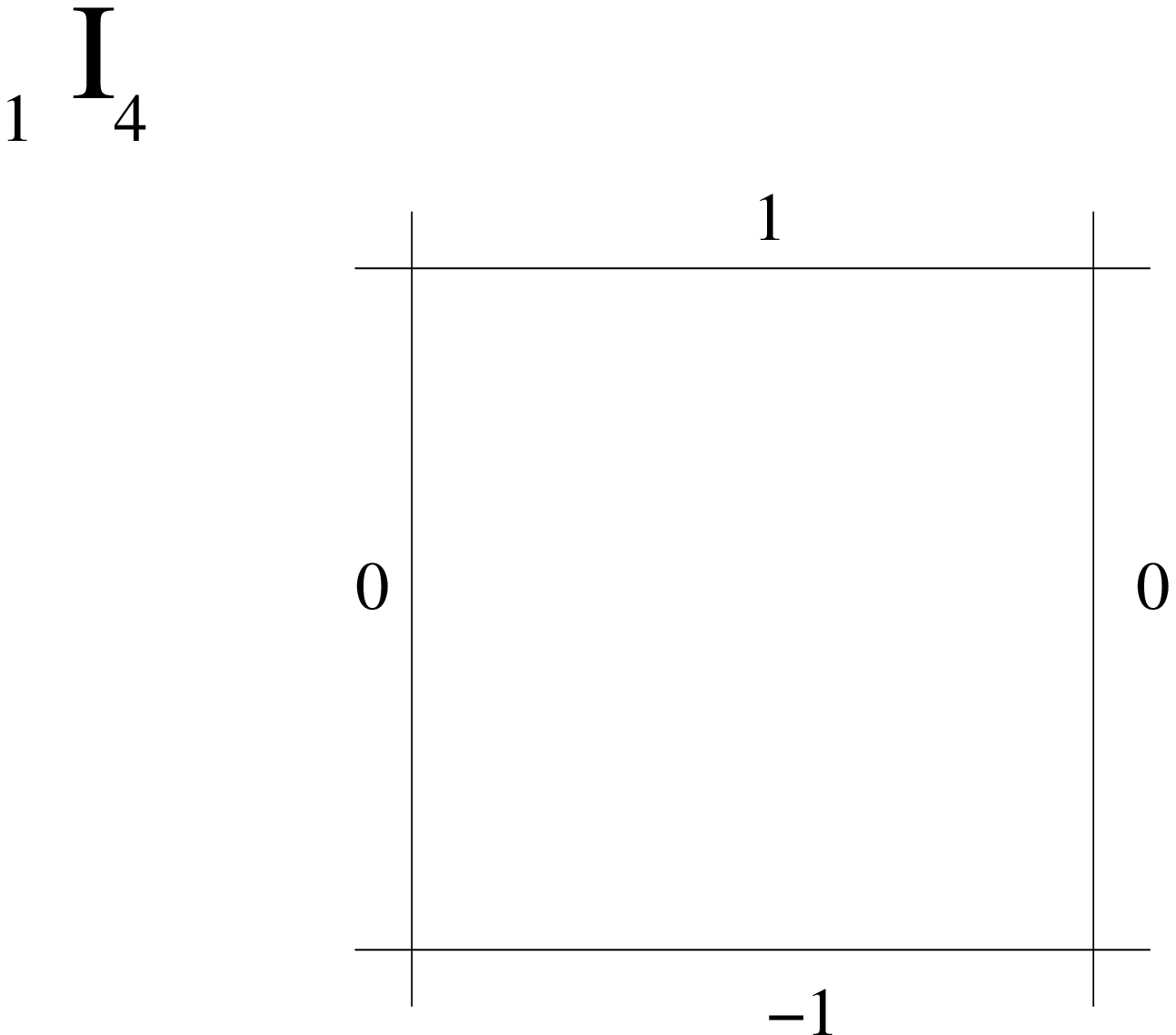}

\includegraphics[scale=0.3]{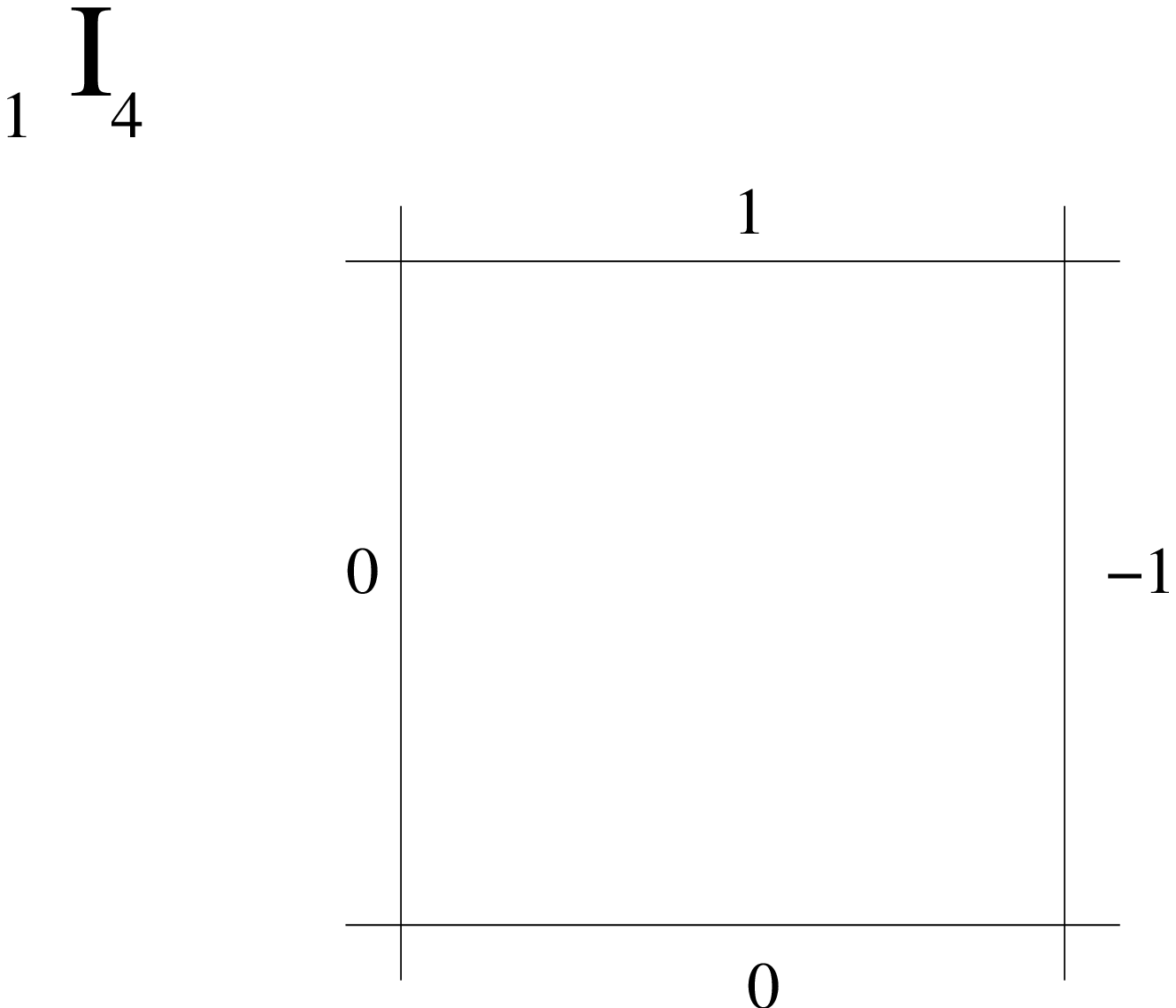}
\hspace{3truecm}
\includegraphics[scale=0.3]{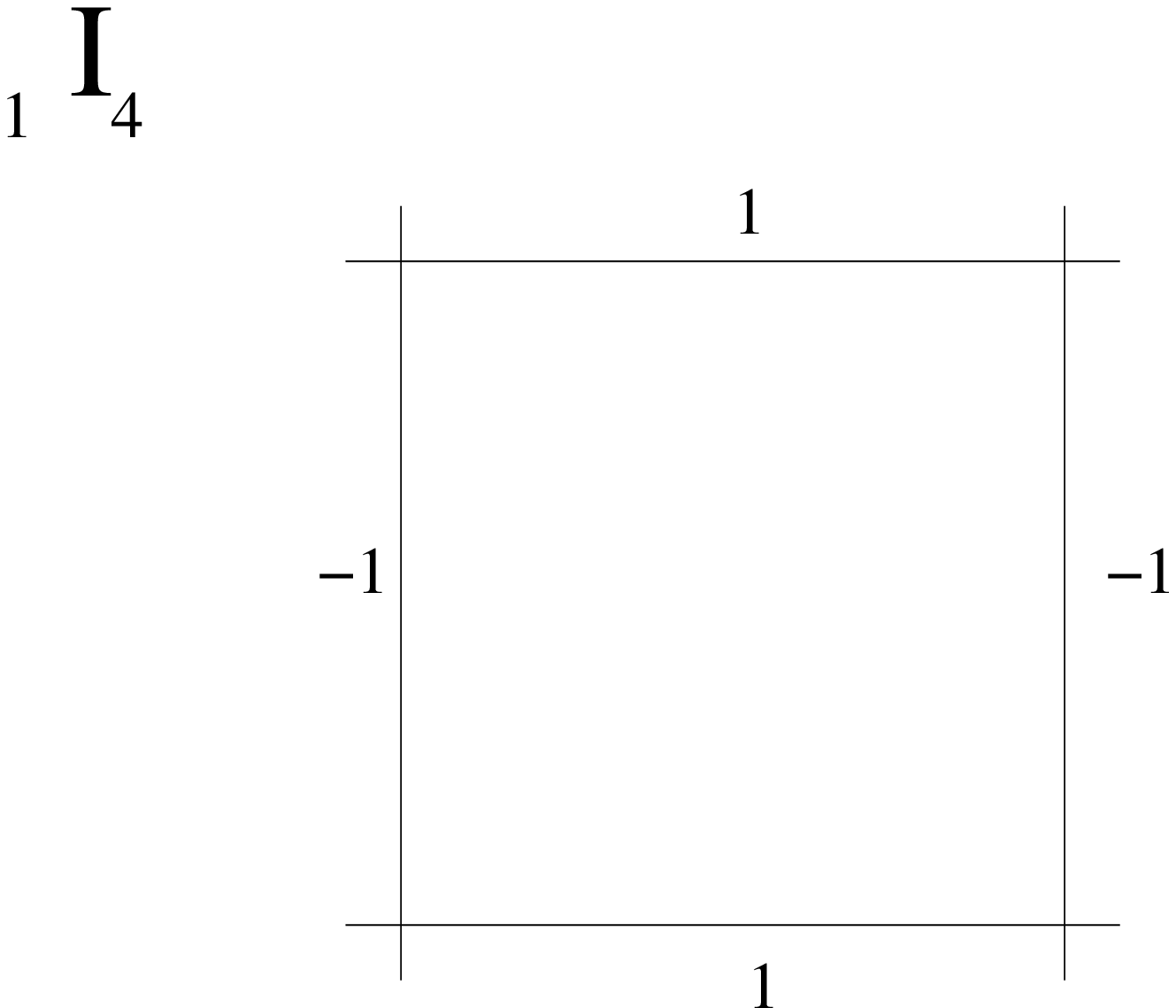}

\noindent{where the numbers denote the degrees of the line bundles
we have to consider on each irreducible component. The following
examples in $I_4$ and $I_6$ are not possible for a strictly
semistable line bundle of degree 0:}}

\

\includegraphics[scale=0.3]{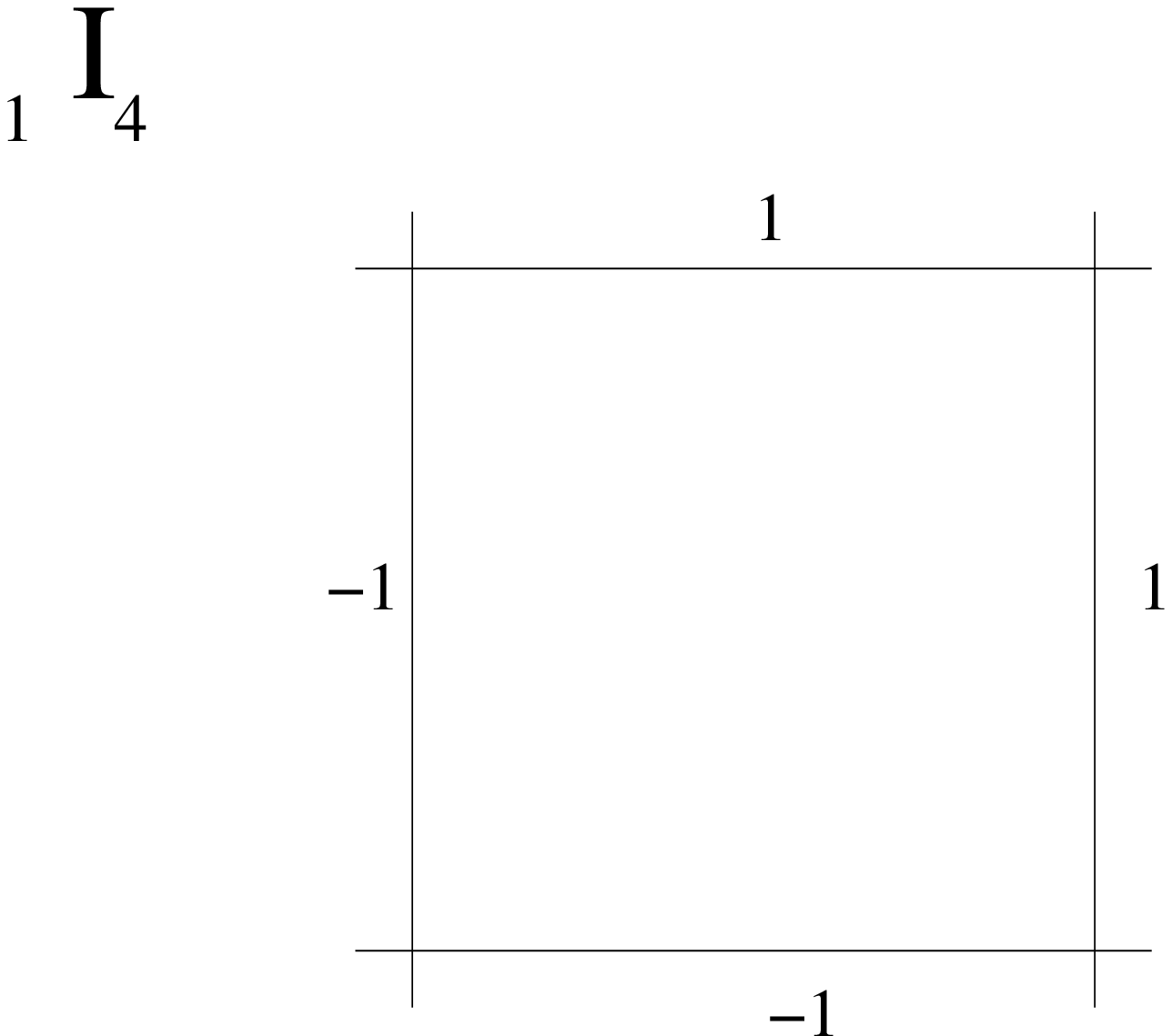}
\hspace{3truecm} \includegraphics[scale=0.3]{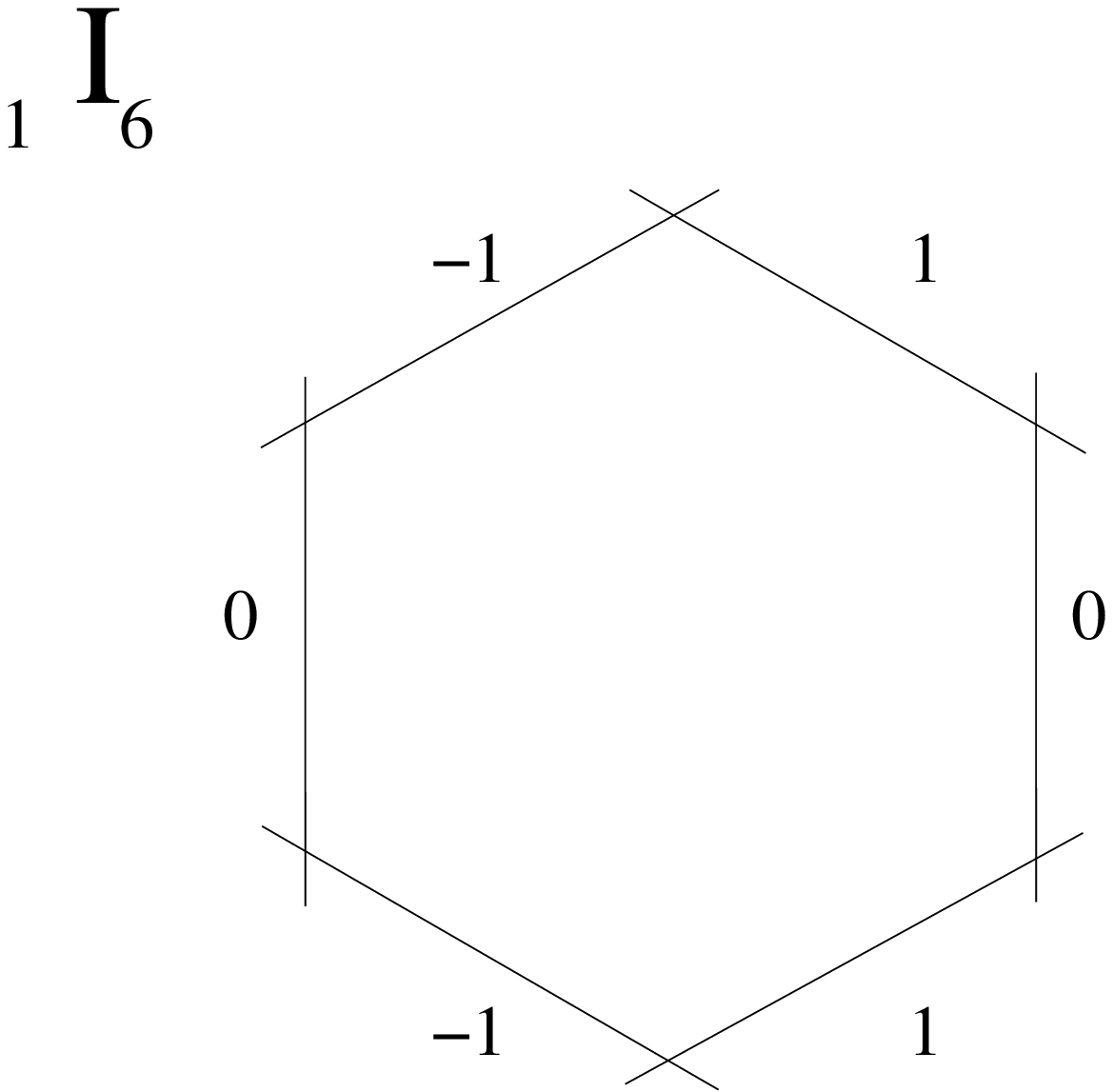}
\end{example}

\begin{corollary} \label{c:nopolarizacion}If $X$ is a polarized curve of type $III$, $IV$
or $I_N$, the (semi)stability  of a pure dimension one sheaf of
rank 1  and degree 0 on $X$ does not depend on the polarization.
\end{corollary}

\begin{proof} For line bundles, the result follows from
Propositions \ref{p:grado0} and \ref{p:propiasemi}. If $F$ is a
pure dimension one sheaf on $X$ of rank 1 and degree 0 which is
not a line bundle, using Lemma \ref{l:buenossubhaces}, it is
semistable if and only if $-\chi(\mathcal{O}_D)\leq d_D$ for any
$D\subset X$ which does not depend on the polarization because
$d_D=\dg_{H_D}(F_D)=\chi(F_D)-\chi(\mathcal{O}_D)$.
\end{proof}

Since we have proved that for these curves there are always
strictly semistable sheaves, we can ensure that

\begin{corollary}\label{c:nofino} Let  $X$ be a curve of type  $III$, $IV$ or $I_N$.
The moduli space of semistable pure dimension one sheaves of rank
1 and degree 0 on $X$ is not a fine moduli space, that is, there
is no universal sheaf.
\end{corollary}

\begin{corollary}\label{c:equivalentes} Let $X$ be a polarized curve of type $III$,
$IV$ or $I_N$. If $F$ is a strictly semistable pure dimension one
sheaf of rank 1 and degree 0 on $X$, then its graded object is
$Gr(F)=\oplus_i\mathcal{O}_{\mathbb{P}^1}(-1)$.
\end{corollary}

\begin{proof} Bearing in mind that for these sheaves it is
$b=h-1$ and $t=-1$, the  result follows from the descriptions
given in the preceding section (see \cite{L} for details).
\end{proof}

\end{document}